\documentclass{article}

\usepackage{arxiv}

\usepackage[utf8]{inputenc} 
\usepackage[T1]{fontenc}    

\usepackage{url}            
\usepackage{booktabs}       
\usepackage{amsfonts}       
\usepackage{amssymb}
\usepackage{amsthm}
\usepackage{nicefrac}       
\usepackage{microtype}      
\usepackage{graphicx}
\usepackage{natbib}

\usepackage{multirow}
\usepackage{booktabs}
\usepackage{graphicx,epstopdf}
\usepackage[caption=false]{subfig}
\usepackage{algorithm}
\usepackage{algorithmic}
\usepackage{float}
\usepackage[table]{xcolor}
\definecolor{regfivetwelve}{RGB}{222,235,247} 
\definecolor{regtenfour}{RGB}{226,240,217}    
\definecolor{regtwentyfour}{RGB}{252,229,214} 
\usepackage{tikz}
\usetikzlibrary{positioning, arrows.meta, shapes, calc, backgrounds, fit}

\usepackage{amsmath,amsfonts,bm}

\def\eqref#1{Eq.~\ref{#1}}

\def\1{\bm{1}}

\DeclareMathAlphabet{\mathsfit}{\encodingdefault}{\sfdefault}{m}{sl}
\SetMathAlphabet{\mathsfit}{bold}{\encodingdefault}{\sfdefault}{bx}{n}

\usepackage{hyperref}       
\usepackage{doi}
\usepackage{cleveref}       
\graphicspath{{../}}

\title{McMg: A Learned Phase-Space Multi-channel Multigrid Preconditioner for Helmholtz Equations}
\renewcommand{\shorttitle}{\textit{arXiv} Template}
\date{}

\author{
  Jiwei Jia\textsuperscript{a,c} \quad
  Xinliang Liu\textsuperscript{b} \quad
  Juntao Wang\textsuperscript{a,c} \quad
  Jinchao Xu\textsuperscript{d}\\[2pt]
  \normalfont\small\textsuperscript{a}\,School of Mathematics, Jilin University, Changchun 130012, China\\
  \normalfont\small\textsuperscript{b}\,Ocean University of China, Qingdao 266100, China\\
  \normalfont\small\textsuperscript{c}\,Shenzhen Loop Area Institute, Shenzhen 518038, China\\
  \normalfont\small\textsuperscript{d}\,King Abdullah University of Science and Technology, Thuwal, Saudi Arabia
}

\hypersetup{
pdftitle={McMg},
pdfauthor={Jiwei Jia, Xinliang Liu, Juntao Wang, Jinchao Xu},
}

\begin{document}
\maketitle

\begin{abstract}
Solving heterogeneous Helmholtz equations at high wavenumbers remains challenging because the discretized operator is indefinite, pollution degrades phase accuracy, and scalar coarse-grid correction can discard the local phase and propagation-direction information carried by oscillatory errors.
We propose Multi-channel Multigrid (McMg), a learned phase-space multigrid preconditioner for heterogeneous Helmholtz equations. Rather than predicting the solution directly, McMg maps residuals to corrections within an iterative framework. Its central idea is to coarsen physical space while retaining unresolved local wave information in the channel dimension: each coarse node carries a learned packet of amplitude, phase, direction, and scattering coefficients rather than a single scalar unknown.
The architecture combines linear multi-channel transfer operators with locally adaptive stencils, neural PDE operators, and medium-dependent smoothers whose coefficients are generated from the wave speed. For a fixed medium, the V-cycle is linear in the residual; nonlinear physical features are computed once in a setup phase and cached, so each online iteration reduces to convolutions with fixed coefficients.
We further study generalization across scales. Models trained on small domains transfer directly to larger domains and higher effective wavenumbers, and a Layer-by-Layer Progressive Finetuning (LLPF) strategy improves large-domain scalability by adding new coarse levels while finetuning only the newly introduced parameters.
Numerical experiments on high-frequency, high-contrast, and large-scale three-dimensional problems demonstrate that McMg requires substantially fewer iterations and less wall-clock time than strong classical baselines, while consistently outperforming existing neural preconditioners.
\end{abstract}

\section{Introduction}
The Helmholtz equation, which governs time-harmonic wave propagation, plays a fundamental role in a wide range of applications, including acoustics, electromagnetics, geophysics, and medical imaging. Despite its importance, computing accurate numerical solutions remains notoriously challenging, particularly in the presence of heterogeneous media and high wavenumbers. At high wavenumbers, standard discretizations suffer from the pollution effect, with a wavenumber-dependent phase error that accumulates across the domain~\cite{babuska1997pollution}.  In multilevel settings this difficulty is amplified: using the same fixed stencil on progressively coarser grids often yields larger dispersion errors, leading to phase mismatches between levels~\cite{stolk2016dispersion}. Meanwhile, the discretized Helmholtz operator is strongly indefinite and may be close to singular in regimes of practical interest, which breaks the assumptions behind classical multigrid. Standard relaxations then cease to be reliable smoothers, and coarse-grid correction may no longer approximate the dominant error components. As a result, multigrid iterations can stagnate or diverge, with the failure often more pronounced in heterogeneous media \cite{ernst2011difficult}.

These challenges have motivated a broad class of Helmholtz-specific solvers and preconditioners. Within the multigrid family, shifted-Laplacian methods~\cite{erlangga2004class} use a complex-shifted operator to add damping, improving spectral properties and yielding effective Krylov preconditioners. Wave-ray multigrid~\cite{brandt1997wave,livshits2006accuracy} augments standard multigrid with ray-based corrections to represent oscillatory components that cannot be resolved on very coarse grids. Beyond multigrid, the Convergent Born Series (CBS)~\cite{osnabrugge2016convergent} constructs a modified Born expansion based on a damped reference medium and a suitable preconditioner, yielding a provably convergent iteration even for large, strongly scattering inhomogeneous media. Domain decomposition methods~\cite{chen2013source,leng2022trace} couple local subdomain solves through inter-subdomain transfer mechanisms, enabling efficient sweeping-type solvers and scalable implementations. Despite these advances, achieving robustness at high wavenumbers and in strongly heterogeneous media without sacrificing efficiency remains an open problem~\cite{ernst2011difficult}.

A useful way to state the multigrid obstruction is that scalar coarsening is not well matched to the phase-space structure of high-frequency waves. After relaxation, a Helmholtz error is not generally smooth in $x$; in the high-frequency regime it is more naturally described by a local WKB/geometrical-optics ansatz $a(x)e^{i\phi(x)}$, with wave vector $\nabla\phi(x)$ near the characteristic shell $|\xi|\approx k(x)$~\cite{babich1991short,engquist2003computational}. Coarsening only the physical coordinate $x$ therefore discards local direction, phase, and scattering information that remains essential for correction. This suggests a different coarse object: a coarse node should carry several local wave-packet or phase-space coefficients, not a single scalar degree of freedom~\cite{ralston1982gaussian}.


Recently, there has been growing interest in leveraging deep learning for solving partial differential equations (PDEs). Physics-Informed Neural Networks (PINNs)~\cite{raissi2019physics} incorporate the governing equations into the training objective and have been applied to a variety of forward and inverse problems, but they typically require retraining for each new parameter configuration and thus generalize poorly across parametrized PDE families. In contrast, neural operator methods~\cite{azizzadenesheli2024neural} aim to learn mappings between infinite-dimensional function spaces. Representative architectures include DeepONet~\cite{lu2021learning}, which combines a branch network encoding the input function and a trunk network encoding spatial coordinates, and the Fourier Neural Operator (FNO)~\cite{li2020fourier}, which performs global spectral convolutions to capture long-range dependencies. Subsequent extensions further broaden applicability, such as Geo-FNO~\cite{li2022fourier} for irregular geometries, and MgNO~\cite{he2024mgno} with a multigrid-inspired hierarchical design; attention-based operator learning has also been explored, e.g., Galerkin-type attention~\cite{cao2021choose}.

Although neural operator methods overcome the generalization limitations of PINNs and enable direct inference, they still face challenges in achieving high accuracy and ensuring solution reliability. This issue stems from the well-known spectral bias of neural networks~\cite{rahaman2019spectral, liu2024mitigating}, which tend to learn low-frequency components efficiently while struggling to capture high-frequency modes. 
To address this, recent studies have proposed Neural Solvers or Preconditioners~\cite{lerer2024multigrid,rudikov2024neural,zhang2024blending,cui2025neural,xiemgcfnn,stanziola2021helmholtz,kopanivcakova2025leveraging,zou2026probabilistic}. Instead of acting as end-to-end predictors, these methods embed neural networks as learned components within classical iterative frameworks, allowing for progressive error reduction.
Notable examples include WAVE-ADR-NS~\cite{cui2025neural}, an efficient multigrid algorithm inspired by the Wave-Ray method~\cite{livshits2006accuracy} that learns to eliminate ``characteristic'' error components specific to the high-wavenumber Helmholtz equation. Similarly, MGCFNN~\cite{xiemgcfnn} constructs a hierarchical AI solver by combining a neural multigrid architecture with a FNO at the coarse level. A closely related family is hybrid methods~\cite{zhang2024blending, lee2025hybrid, lee2025fast}, which interleave classical iterative solvers with neural operators, using the network to capture the low-frequency error that the conventional smoother struggles to eliminate and thereby accelerate its convergence.
While promising, many neural solvers for Helmholtz equations still face practical limitations. Although they reduce the iteration count relative to classical methods, their accuracy and convergence often degrade at high wavenumbers and in strongly heterogeneous media. Moreover, the computational cost of evaluating relatively heavy networks at each iteration can offset these gains, leaving the overall solving time as slow as, or even slower than, that of classical methods.

In this work, we develop Multi-channel Multigrid (McMg), a learned phase-space multigrid preconditioner for heterogeneous Helmholtz equation. The method coarsens physical space as in a standard multigrid hierarchy, but keeps a finite-dimensional representation of unresolved phase, direction, and scattering information in the channel dimension. Our objective is to address a key limitation of existing neural preconditioners: whether efficient solution performance can be retained for high-frequency, high-contrast, and large-scale Helmholtz problems, as required in realistic industrial applications.
Our main contributions are summarized as follows:
\begin{itemize}
    \item We introduce McMg, a learned phase-space multi-channel multigrid preconditioner for heterogeneous Helmholtz problems. The method is formulated as an iterative error-correction scheme rather than an end-to-end neural predictor, and learns an approximation to the inverse operator that can be used as a solver or preconditioner.

    \item We propose multi-channel multigrid operators based on locally adaptive stencils. The stencil coefficients are generated from the local wave speed, while the residuals and corrections are lifted to a multi-channel space. The channel dimension acts as a learned local packet basis, helping preserve phase and directional information that scalar coarse grids lose.

    \item We develop a setup strategy with parameter sharing within the multigrid hierarchy. Since the medium coefficient is fixed during the solve, the nonlinear physical features are computed once in the setup phase and cached; every subsequent iteration then applies the cached adaptive operators to the evolving residuals and corrections as a sequence of linear convolutions. Thus, for a fixed medium, the learned V-cycle is linear in the residual.

    \item We demonstrate cross-scale generalization and large-domain scalability. Models trained on small domains or patches can be deployed on substantially larger domains and higher effective wavenumbers. For further scaling, we introduce Layer-by-Layer Progressive Finetuning (LLPF): by freezing previously trained levels and finetuning only a single newly appended coarse level, LLPF extends the long-range support of the solver to large domains at a fraction of the cost of full retraining, while matching its convergence.

    \item We validate McMg on high-frequency, high-contrast, and large-scale Helmholtz problems, spanning the OpenBreastUS, OpenFWI CurveFault-B, and Kimberlina 3D benchmarks. Relative to the fastest classical baseline on each benchmark---selected among GMRES with shifted-Laplacian multigrid, CBS, and sparse LU---McMg reduces iteration counts by $110.6\times$, $51.3\times$, and $45.1\times$, with corresponding wall-clock speedups of $10.3\times$, $3.4\times$, and $6.8\times$. Spectral-error, coarse-basis, channel-SVD, and Green's-function-support diagnostics provide evidence for the learned phase-space interpretation.
\end{itemize}

\section{Preliminaries}
In this section, we establish the mathematical formulation of the heterogeneous Helmholtz equation, discuss the absorbing boundary layers required for domain truncation, and review its finite difference discretization. We then outline the V-cycle Geometric Multigrid algorithm.
\subsection{Heterogeneous Helmholtz equation}
We consider the Helmholtz equation, which models wave scattering phenomena in an unbounded heterogeneous medium. To ensure the well-posedness of the solution, we assume that the equation satisfies the Sommerfeld radiation condition. The complete formulation is given by:
\begin{equation}\label{eq:helm}
\begin{aligned}
    -\Delta u(x) - k(x)^2 u(x) &= f(x), \quad x \in \mathbb{R}^{D},\\
    \lim_{|x| \to \infty} |x|^{\frac{D-1}{2}} \left( \frac{\partial u}{\partial |x|} - i k_0 u \right) &= 0,
\end{aligned}
\end{equation}
where $u(x)$ denotes the wavefield, $f(x)$ is the source term, $D$ represents the spatial dimension, and $i = \sqrt{-1}$ is the imaginary unit. The spatially varying wavenumber is given by $k(x) = \frac{\omega}{c(x)}$, where $\omega$ is the angular frequency and $c(x)$ is the speed of sound. In the radiation condition, $k_0$ denotes the constant background wavenumber, assuming $k(x) \to k_0$ as $|x| \to \infty$.

To facilitate numerical computation, it is necessary to truncate the unbounded domain. However, applying Dirichlet or Neumann boundary conditions on the artificial boundary typically leads to spurious reflections, which do not align with the physical behavior of wave propagation in unbounded domains. To mitigate such reflections, absorbing boundary layers are employed; common techniques include the sponge layer~\cite{israeli1981approximation} and the perfectly matched layer (PML)~\cite{berenger1994perfectly,stanziola2021helmholtz}. 

The sponge layer adds a damping term to the original equation,
\begin{equation}
    -\Delta u(x) - k(x)^2 (1-\gamma i)u(x) = f(x), \quad x \in \Omega \cup \Omega_{\mathrm{abl}},
\end{equation}
where $\Omega$ denotes the computational domain and $\Omega_{\mathrm{abl}}$ the absorbing boundary layer; the damping coefficient $\gamma$ vanishes inside $\Omega$ and increases polynomially with distance from the interface within $\Omega_{\mathrm{abl}}$. The PML instead extends the spatial coordinates into the complex plane, introducing artificial decay without reflection through a complex stretching of the Laplacian; we refer the reader to Appendix~\ref{sec:pml_derivation} for the detailed formulation.

In this work, the computational domain is chosen as a rectangular box (a cube in the three-dimensional case) and discretized on a structured uniform grid of spacing $h$. We adopt the finite difference method (FDM) as the underlying discretization; although our method ultimately learns adaptive stencils, it builds upon this classical scheme. Due to the heterogeneous wavenumber $k(x)$, the resulting discrete operator acts as a spatially varying convolution:
\begin{equation}\label{eq:discrete_fdm}
    \mathcal{K}_{A} * \bm u = \bm f,
\end{equation}
where $\mathcal{K}_{A}$ denotes the local stencil, $*$ the spatially variant convolution, $\bm u$ the numerical solution, and $\bm f$ the discrete source term. In two dimensions, the standard five-point stencil at grid index $(i,j)$ takes the form:
\begin{equation}\label{eq:FDM_stencil}
    \mathcal{K}_A[i,j] = \frac{1}{h^2}
    \begin{bmatrix}
        0 & -1 & 0 \\
        -1 & 4 - (\bm k_{i,j} h)^2 & -1 \\
        0 & -1 & 0
    \end{bmatrix},
\end{equation}
where $\bm k_{i,j}$ denotes the discrete wavenumber $\bm k$ at location $(i,j)$. Other discretizations, such as the Fourier spectral method (Appendix~\ref{sec:app_fsm}), are equally applicable. Regardless of the scheme, the resulting linear system can be written in the general form
\begin{equation}
    A \bm u = \bm f.
\end{equation}

\subsection{Geometric Multigrid (GMG)}
Geometric multigrid (GMG) is an optimal-complexity ($\mathcal{O}(N)$) iterative solver for elliptic PDEs that combines fine-grid smoothing with coarse-grid correction over a hierarchy of nested grids~\cite{xu1997introduction}. Its efficiency stems from a complementary division of labor: inexpensive local smoothers (such as weighted Jacobi or Gauss--Seidel) rapidly damp the high-frequency components of the error, while the remaining smooth, low-frequency error is transferred to a coarser grid, where it is cheaper to resolve and, relative to the coarser mesh, once again appears oscillatory. Applying this principle recursively yields the V-cycle: starting from the finest grid, the residual is smoothed and then restricted down the hierarchy to the coarsest level $l=L$, where the system is solved directly; the resulting correction is prolongated back up, with smoothing applied at each level along the way.

While highly effective for elliptic operators, standard GMG can deteriorate or diverge for high-frequency Helmholtz problems~\cite{ernst2011difficult}. Applying the same fixed stencil across levels yields inconsistent dispersion on coarser grids, so coarse-grid corrections may arrive with an incorrect phase and amplify the error instead of reducing it.

\section{Multi-channel Multigrid}

In this section, we detail the proposed Multi-channel Multigrid (McMg) architecture. Unlike classical geometric multigrid methods that operate on scalar fields across levels, McMg uses a multi-channel latent representation. The spatial grid is coarsened, but each coarse node carries several latent coefficients that can encode local phase, propagation direction, and scattering information. Throughout the hierarchy, the discrete operators and smoothers are modulated by the local wave speed, so that the solver adapts to heterogeneous media while preserving wave information often lost in standard scalar GMG.

\subsection{Learnable Setup Phase}
The operators used by McMg share a common structure. Given the medium coefficient $\bm{m}$ (e.g., the discrete wavenumber $\bm{k}$), a learnable nonlinear map first extracts a physical feature
\begin{equation}\label{eq:phys_feature}
    \bm{c} = \sigma(\mathcal{K}_{\bm m} \ast \bm{m}),
\end{equation}
which then modulates a linear convolution of the field $\bm{x}$ (a residual or correction):
\begin{equation}\label{eq:pyhs_conv}
    \Phi_{\bm m}(\bm{x}) = \bm{c} \odot (\mathcal{K}_{\bm x} \ast \bm{x}).
\end{equation}
Here $\mathcal{K}_{\bm m}$ and $\mathcal{K}_{\bm x}$ are learnable multi-channel convolution kernels, $\ast$ denotes multi-channel convolution, $\odot$ is element-wise multiplication, and $\sigma$ is a nonlinear activation (GELU~\cite{hendrycks2016gaussian} in this work). The construction is deliberately asymmetric in its two inputs. It is \emph{nonlinear} in the coefficient $\bm{m}$, so that the physical feature $\bm{c}$ can represent the complex dependence of optimal discretization parameters on the wavenumber (e.g., the dispersion-optimized stencil coefficient $2\cos(\bm{k}h)$, which is transcendental in $\bm{k}$ rather than a low-order polynomial). For a fixed medium, however, $\Phi_{\bm m}$ is \emph{linear} in the field $\bm{x}$, matching the linearity of the Helmholtz solution operator $A^{-1}$: the feature $\bm{c}$ acts as a spatially adaptive stencil whose weights are set by the medium and applied linearly to $\bm{x}$.

A key consequence of this split is that the physical feature $\bm{c}$ in \eqref{eq:phys_feature} depends only on the medium coefficient, which is fixed for a given problem. It can therefore be computed once and reused at every iteration of the V-cycle, mirroring how a classical multigrid method assembles its operator hierarchy once before iterating. The learnable setup phase precomputes and caches these features for all levels, so that each subsequent iteration reduces to a sequence of linear convolutions against the cached coefficients, as in \eqref{eq:pyhs_conv}. This keeps the heavy, nonlinear feature extraction out of the iterative loop and substantially lowers the per-iteration cost. For brevity, we omit the parameter subscript $\theta$ in the learnable operators below.

\paragraph{Feature Hierarchy Construction} We generate two distinct sets of cached physical features for each grid level $l$:
\begin{enumerate}
    \item \textbf{PDE Features ($\bm a_l$):} These feature maps act as a learned, high-dimensional generalization of finite difference stencil coefficients. They encode local dispersion corrections and wave speeds required to approximate the Helmholtz operator $A_l$.
    \item \textbf{Smoother Features {($\bm s_l$)}:} These maps serve as spatially varying relaxation parameters. They allow the solver to dynamically adjust smoothing intensity based on local heterogeneity.
\end{enumerate}

The detailed construction of the multigrid hierarchy is summarized in Algorithm~\ref{alg:McMg_setup}. To facilitate the generation of physical features, we employ a set of non-linear convolutional neural networks (CNNs) at each level. Specifically, let $f_{\text{lift}}$ denote the lifting network that maps the raw physical coefficients into the high dimensional latent feature space. For level $l$, we define $f_{\bm s}^{l}$ as the network projecting PDE features $\bm a_l$ into smoother features $\bm s_l$, and $f_{\bm a}^{l}$ as the restriction network mapping features $a_l$ to the coarser level. All three types of networks ($f_{\text{lift}},f_{\bm s}^{l},f_{\bm a}^{l}$) are composed of multi-channel non-linear convolutional layers. Conversely, the inter-grid transfer operators $\mathcal R_{l}^{l+1}$ and $\mathcal P_{l+1}^{l}$ are strictly multi-channel linear convolutional (or transpose convolutional) layers to preserve the linearity of the solver.
\begin{algorithm}[htbp]
\caption{McMg Learnable Setup}\label{alg:McMg_setup}
\begin{algorithmic}[1]
\REQUIRE Fine grid discrete wavenumber $\bm k$, maximum level $L$.
\ENSURE Hierarchy of physical features $\{\bm a_l, \bm s_l\}_{l=1}^L$
\STATE \textit{// Apply lifting}
\STATE $\bm a_1\leftarrow f_{\text{lift}}(\bm k)$
\STATE \textbf{Hierarchy Construction}
\FOR{$l = 1$ to $L$}
    \STATE \textit{// Generate smoother features} 
    \STATE $\bm s_l \leftarrow f_{\bm s}^{l}(\bm a_l)$
    \IF{$l < L$}
        \STATE \textit{// Generate coarse level PDE features}
        \STATE $\bm a_{l+1} \leftarrow f_{\bm a}^{l}(\bm a_l)$
    \ENDIF
\ENDFOR
\RETURN $\{\bm a_l, \bm s_l\}_{l=1}^L$
\end{algorithmic}
\end{algorithm}

\subsection{Learnable Solving Phase}
Building upon the standard Geometric Multigrid framework~\cite{xu1997introduction}, the solving phase of McMg executes a recursive V-cycle comprising three core operations: neural smoothing, neural restriction, and neural prolongation. These operations reuse the cached physical features precomputed during the learnable setup phase, modulating linear convolutions of the residuals and corrections as in \eqref{eq:pyhs_conv}. We denote the complete neural V-cycle as $\mathcal{MG}(\bm{k}, \bm{r})$, where $\bm{k}$ represents the discrete wavenumber, and $\bm r$ is the input residual. Since the physical coefficients remain fixed throughout the iterative process, we simplify the notation to $\mathcal{MG}(\bm{r})$ for brevity. 

The overarching objective of the McMg is to function as a high-quality preconditioner that approximates the inverse Helmholtz operator:
\begin{equation}
    \mathcal{MG}(\bm{r}) \approx A^{-1}  \bm r.
\end{equation}
Crucially, since the inverse operator $A^{-1}$ is linear, McMg is designed to preserve linearity with respect to $\bm r$ for fixed $\bm k$. Nonlinearity enters only through the medium-dependent setup features, not through the residual-to-correction map used during the solve.

\paragraph{Neural Smoothing}
In contrast to fixed relaxation schemes like Jacobi or Gauss-Seidel, McMg employs a learnable neural smoother built from the cached smoother features $\bm s_l$. At grid level $l$ and smoothing iteration $i$, the neural smoother is defined as:
\begin{equation}
    \mathcal{S}_{l}^{(i)}(\bm x) = \mathcal{W}^{(i)} \ast \big(\bm s_l \odot (\mathcal{K}_{\bm x} \ast \bm x)\big),
\end{equation}
where $\mathcal{W}^{(i)}$ denotes the kernel of a standard multi-channel linear convolution.
Numerical experiments demonstrate that employing iteration-dependent smoother parameters significantly accelerates convergence. Consequently, while the smoother features are shared, we assign a unique set of weights $\mathcal{W}^{(i)}$ for each smoothing step $i$. Conversely, to ensure physical consistency across iterations, the wave propagation physics are governed by a single Neural PDE Operator $\mathcal{A}_{l}$ shared across the entire level:
\begin{equation}
\mathcal{A}_{l}(\bm{x}) = \mathcal{W}_l \ast \bm{x} + \bm{a}_l \odot (\mathcal{K}_{\bm x} \ast \bm{x}),
\end{equation}
where the first term employs a linear convolution with kernel $\mathcal{W}_l$ to approximate the Laplacian operator, and the second term accounts for the heterogeneous physics via the precomputed PDE features $\bm{a}_l$.
Combining these operators, the neural smoothing phase is updated via:
\begin{equation}
\bm{e}_l^{(i)} = \bm{e}_l^{(i-1)} + \mathcal{S}_{l}^{(i)} \left( \bm{r}_l - \mathcal{A}_{l}(\bm{e}_l^{(i-1)}) \right).
\end{equation}
Crucially, these operations are performed within a high-dimensional feature space. Unlike the scalar fields used in classical methods, the multi-channel features $(\bm{e}_l,\bm{r}_l)$ attach several values to each grid point. The smoother therefore acts not only as a damping step, but also as a medium-dependent channel-mixing operator that can locally modify amplitude, phase, and direction of latent wave packets.

\paragraph{Neural Restriction}
The objective of the restriction step is to transfer residuals from the fine grid to the coarse grid while preserving the components that are difficult to eliminate locally. While classical geometric multigrid relies on fixed scalar kernels, our neural restriction generalizes this operation to a learnable linear operator implemented as a multi-channel strided convolution, denoted by $\mathcal{R}_{l}^{l+1}$. This allows the network to project fine-grid latent residual features onto a coarse latent space that can still distinguish phase and directional components. The coarse-level latent residual features $\bm r_{l+1}$ are computed as:
\begin{equation}
    \bm r_{l+1} = \mathcal R_{l}^{l+1} \left( \bm r_l - \mathcal A_{l}(\bm e_l) \right),
\end{equation}


\paragraph{Neural Prolongation}
Analogous to the neural restriction, the neural prolongation step employs a learnable linear operator implemented as a multi-channel transposed convolution, denoted by $\mathcal{P}_{l+1}^{l}$. This operator maps coarse-grid latent error-correction features back to the fine-grid feature space. Unlike standard fixed interpolation, the learned prolongation defines a multi-channel packet dictionary: one coarse spatial node corresponds to several fine-grid oscillatory responses. The fine-grid latent correction features are updated via:
\begin{equation}
    \bm e_{l} \leftarrow \bm e_{l} + \mathcal P_{l+1}^{l} \left( \bm e_{l+1}\right).
\end{equation}

\begin{figure}
    \centering
    \includegraphics[width=\linewidth]{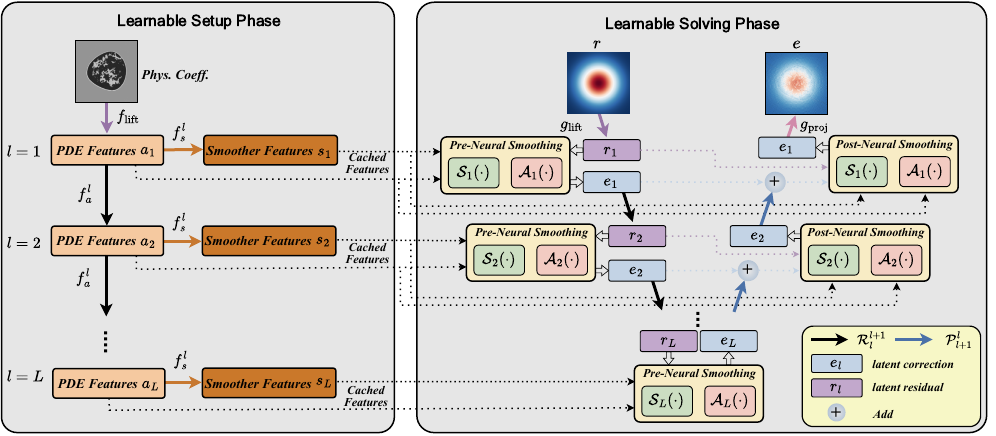}
    \caption{McMg architecture. The setup phase precomputes and caches per-level multi-channel physical features, while the solving phase reuses them within the learnable V-cycle to compute each iterative correction efficiently.}
    \label{fig:mpmg_arch}
\end{figure}

Integrating these learnable components, we construct the complete McMg V-cycle. This framework propagates multi-channel phase-space information across the grid hierarchy rather than relying on scalar coarse-grid variables. The overall solving procedure is summarized in Algorithm~\ref{alg:McMg_solve}. Note that the algorithm employs two linear convolutional modules to interface with the physical field: an encoder $g_{\text{lift}}$ that maps the initial scalar residual into the high-dimensional latent space, and a decoder $g_{\text{proj}}$ that reconstructs the final scalar error correction from the latent features. These lifting and projection maps contain no nonlinear activations, so they preserve linearity with respect to the residual input. The overall procedure is illustrated in Fig.~\ref{fig:mpmg_arch}.
\begin{algorithm}
\caption{McMg Learnable V-Cycle Solving}
\label{alg:McMg_solve}
\begin{algorithmic}[1]
\REQUIRE Current residual $\bm r$, maximum level $L$, cached physical features $\{\bm a_l ,\bm s_l\}_{l=1}^L$, initialized linear transfer operators $\{\mathcal R_{l}^{l+1}, \mathcal P_{l+1}^l\}_{l=1}^{L-1}$, pre-smoothing steps $\{\nu_1^l\}_{l=1}^L$, post-smoothing steps $\{\nu_2^l\}_{l=1}^L$.
\ENSURE Error correction $\bm e =\mathcal{MG}(\bm r)\approx A^{-1}\cdot \bm r$.
\STATE \textit{// Apply lifting}
\STATE $\bm r_1 = g_{\text{lift}}(\bm r)$
\STATE \textbf{Down Cycle:}
\FOR{$l = 1$ to $L$}
    \STATE $\bm e_l^{(0)} = \bm 0$.
    \STATE \textit{// Pre-Neural Smoothing}
    \FOR{$i = 1$ to $\nu_1^l$}
        \STATE $\bm e_l^{(i)} = \bm e_l^{(i-1)} + \mathcal S_{l}^{(i)} \left( \bm r_l - \mathcal A_{l}\left(\bm e_l^{(i-1)}\right) \right)$
    \ENDFOR
    \STATE $\bm e_l = \bm e_l^{(\nu_1^l)}$

    \IF{$l < L$}
        \STATE \textit{// Neural Restriction} 
        \STATE $\bm r_{l+1} = \mathcal R_{l}^{l+1} \left( \bm r_l - \mathcal A_{l}\left(\bm e_l\right) \right)$
    \ENDIF
\ENDFOR

\STATE \textbf{Up Cycle:}
\FOR{$l = L-1$ to $1$}
    \STATE \textit{// Neural Prolongation} 
    \STATE $\bm e_{l} =\bm e_{l}+ \mathcal P_{l+1}^{l} \left( \bm e_{l+1}\right)$
    \STATE \textit{// Post-Neural Smoothing}
    \STATE $\bm e_l^{(0)} = \bm e_l$
    \FOR{$i = 1$ to $\nu_2^l$}
        \STATE $\bm e_l^{(i)} = \bm e_l^{(i-1)} + \mathcal S_{l}^{(\nu_1^l+i)} \left( \bm r_l - \mathcal A_{l}\left(\bm e_l^{(i-1)}\right) \right)$
    \ENDFOR
    \STATE $\bm e_l = \bm e_l^{(\nu_2^l)}$
\ENDFOR
\STATE \textit{// Apply Projection} 
\STATE $\bm e = g_{\text{proj}}(\bm e_1)$
\RETURN $\bm e$
\end{algorithmic}
\end{algorithm}

\subsection{Training and Inference Strategies}\label{sec:mpmg_iteration_scheme_all}\label{sec:mpmg_iteration_scheme}
McMg is trained by minimizing the residual of the discrete Helmholtz system, which we describe below. The same learned operator can also be embedded within the Born series, yielding the Born+McMg variant which follows the neural preconditioning approach of NPBS~\cite{wang2026neural}, we defer its derivation to Appendix~\ref{sec:app_born_series}.

We consider the discrete Helmholtz system
\begin{equation}
    A\bm u = \bm f,
\end{equation}
where $A$ denotes the discretized Helmholtz operator and $\bm f$ is the right-hand side. Given an approximate solution $\bm u^{(k)}$ with residual $\bm r_{\mathrm{H}}^{(k)}$, the most direct use of McMg is as a stationary correction scheme that updates the solution directly by the neural correction at each iteration:
\begin{equation}\label{eq:mcmg_iteration}
    \bm r_{\mathrm{H}}^{(k)} = \bm f - A\bm u^{(k)},
    \qquad
    \bm u^{(k+1)}
    =
    \bm u^{(k)}
    +
    \mathcal{MG}\big(\bm r_{\mathrm{H}}^{(k)}\big).
\end{equation}
In this setting, McMg itself defines the iterative solver.
To train McMg, following \cite{xiemgcfnn,wang2026neural}, we sample random residuals
\begin{equation}
    \bm r \sim \mathcal{N}(\bm 0,\bm I),
\end{equation}
where $\mathcal{N}$ denotes the multivariate standard normal distribution, and minimize the relative residual of the predicted correction:
\begin{equation}\label{eq:loss_fun}
    \mathcal{L}_{\mathrm{H}}
    =
    \frac{1}{B}
    \sum_{i=1}^{B}
    \frac{
    \left\|
    A_{(i)}
    \mathcal{MG}(\bm k_{(i)}, \bm r)
    -
    \bm r
    \right\|_2
    }{
    \left\|
    \bm r
    \right\|_2
    },
\end{equation}
where $B$ is the batch size and $i$ denotes the sample index. Here, $A_{(i)}$ is the discrete operator for the $i$-th sample, parameterized by the wavenumber $\bm k_{(i)}$. We retain the explicit dependence $\mathcal{MG}(\bm k_{(i)}, \bm r)$ rather than the simplified form $\mathcal{MG}(\bm r)$, since each sample is associated with a distinct wavenumber. The residuals $\bm r$ are resampled at each training iteration to expose the model to a broad distribution of residual modes. This objective trains McMg to approximate the inverse of the Helmholtz operator.
The same learned operator can also serve as a preconditioner within GMRES.

\subsection{The Multi-channel Coarse Space}\label{sec:learned_coarse_space}
The defining feature of McMg, relative to classical scalar multigrid, is that residuals, corrections, and the coarse spaces themselves all reside in a multi-channel latent space rather than on scalar grids. Although wide feature maps are routine in neural networks, here the extra channels serve a specific numerical purpose: they allow a coarse grid to retain local phase-space information that scalar coarsening is forced to discard.

The local phase-space motivation is as follows. In the high-frequency regime, Helmholtz fields and errors are often locally approximated by WKB/geometrical-optics ansatz functions~\cite{babich1991short,engquist2003computational}
\[
    e(x)\approx a(x)e^{i\phi(x)},\qquad |\nabla\phi(x)|\approx k(x).
\]
Thus the unresolved state includes not only the position $x$ but also the local wave vector $\xi=\nabla\phi(x)$. A scalar coarse node cannot distinguish wave packets traveling in different directions or with different phases after the spatial grid has been coarsened; this phase-space viewpoint is also consistent with Gaussian-beam and wave-packet descriptions of high-frequency wave propagation~\cite{ralston1982gaussian,engquist2003computational}. In McMg, the channel index provides a finite-dimensional local representation of these unresolved variables.

Let $\alpha=(\mathbf{i},c)$ denote a degree of freedom at the coarsest level $L$, where $\mathbf{i}$ is a coarse-grid location and $c$ is a channel index. If $\bm\delta_{\alpha}$ is the corresponding one-hot tensor, the raw transfer-defined fine-grid response is obtained by the composite learned prolongation
\begin{equation}\label{eq:coarse_basis}
    \Psi_{\alpha}^{\mathrm{raw}} = g_{\text{proj}}\left(\mathcal{P}_{2}^{1} \circ \mathcal{P}_{3}^{2} \circ \cdots \circ \mathcal{P}_{L}^{L-1} (\bm \delta_{\alpha})\right),
\end{equation}
where $g_{\text{proj}}$ maps the latent correction back to the physical scalar field. This raw basis furnishes a learned tentative packet dictionary shared across media. The basis effectively used by a V-cycle is further modified by the medium-dependent smoother and neural PDE operators. If $\mathcal H_l(\bm k)$ denotes the homogeneous linear part of the level-$l$ update induced by the cached features $\bm a_l,\bm s_l$, then a schematic effective response is
\begin{equation}\label{eq:effective_packet_basis}
    \Psi_{\alpha}^{\mathrm{eff}}(\bm k)
    \approx
    g_{\text{proj}}\,
    \mathcal H_1(\bm k)\mathcal P_2^1
    \mathcal H_2(\bm k)\mathcal P_3^2
    \cdots
    \mathcal P_L^{L-1}\bm\delta_{\alpha}.
\end{equation}
Thus the transfer operators provide a learned multi-channel packet dictionary, while the cached medium-dependent smoothers and neural PDE operators transform it into an operator-adapted effective coarse representation. This is analogous in spirit to smoothed aggregation AMG, except that the tentative basis is multi-channel and oscillatory rather than scalar and smooth.

For classical piecewise-polynomial bases (e.g., $P_1$ elements in standard GMG), the coarse mesh size $H$ is constrained by scalar sampling and dispersion requirements; suppressing the accumulated phase error (the pollution effect) can demand conditions such as $H k^2 \lesssim C_0$~\cite{babuska1997pollution}. Once $H$ is too large relative to the wavelength, a basis with a single scalar degree of freedom per node can no longer approximate oscillatory error components, and coarse-grid correction becomes ineffective. The standard remedy is to enrich the coarse space: advanced multiscale methods such as Localized Orthogonal Decomposition (LOD)~\cite{gallistl2015stable,peterseim2017eliminating,hauck2021multi,freese2024super,lu2025two} and the Generalized Multiscale Finite Element Method (GMsFEM)~\cite{chung2014generalized,ma2023wavenumber,jin2024robust} attach several problem-dependent oscillatory basis functions to each coarse node, which improves the resolution condition.

\begin{figure}[htbp]
  \centering
  \includegraphics[width=1\linewidth]{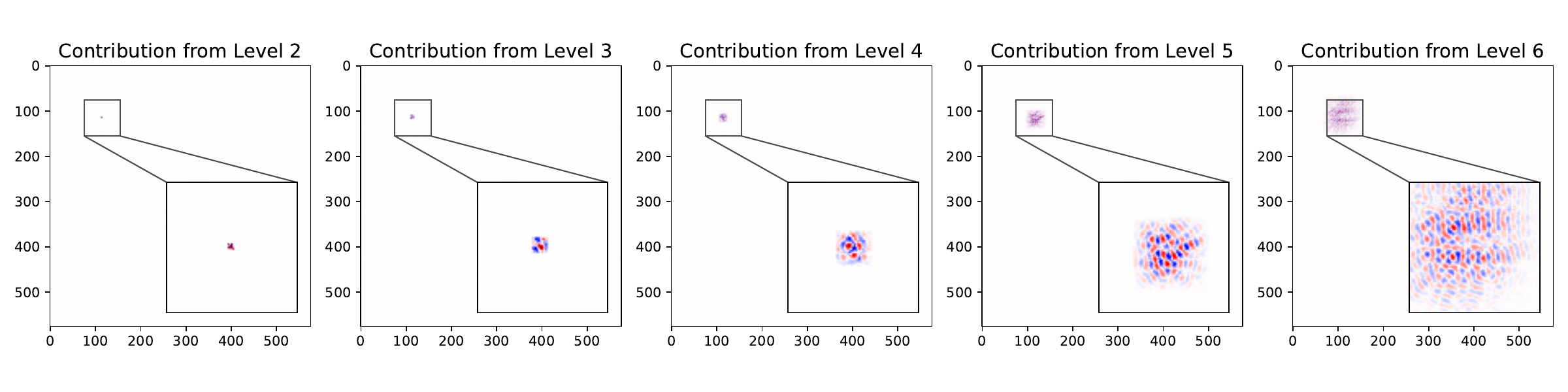}
  \caption{Visualization of the learned coarse-space basis responses. We inject a unit impulse at different coarse levels $l$ and propagate it to the fine grid via \eqref{eq:coarse_basis}. Unlike standard bilinear interpolation, the learned responses (especially at level 6) exhibit rich oscillatory structures spanning many wavelengths. Because each coarse node carries $C$ channels, it contributes a multi-function local packet basis rather than a single scalar interpolation weight; the coarse spatial mesh can be much larger than the wavelength because unresolved phase information is represented in the channel dimension.}
  \label{fig:exp_levels_basis}
\end{figure}

McMg follows the same enrichment principle, with the per-node basis learned end-to-end rather than constructed analytically. We visualize the learned bases for a six-level model trained on the OpenBreastUS dataset under the configuration of Section~\ref{sec:comp_with_other_model_OpenBreastUS}, with one coarse level added beyond the five-level solver used there so as to expose the behavior of the coarsest scale; here $k_0 D \approx 500$ and the coarsest level $L=6$ has spacing $H \approx 32h \approx 5\lambda$. A scalar piecewise-polynomial coarse space would be far too sparse in this regime. As shown in Figure~\ref{fig:exp_levels_basis}, activating a single coarse degree of freedom at level $L=6$---a unit impulse in one channel $c$ at one coarse location $\mathbf{i}$---induces a fine-grid response that is not a localized interpolation kernel but a highly oscillatory wave packet spanning roughly $12$ wavelengths. The coarse representation stays expressive not by reducing $H$, but by enriching each node with $C$ learned packet coefficients.


\begin{figure}[htbp]
  \centering
  \subfloat[Per-channel fine-grid responses.\label{fig:similarity_channel_basis_response}]{%
    \includegraphics[width=0.7\linewidth]{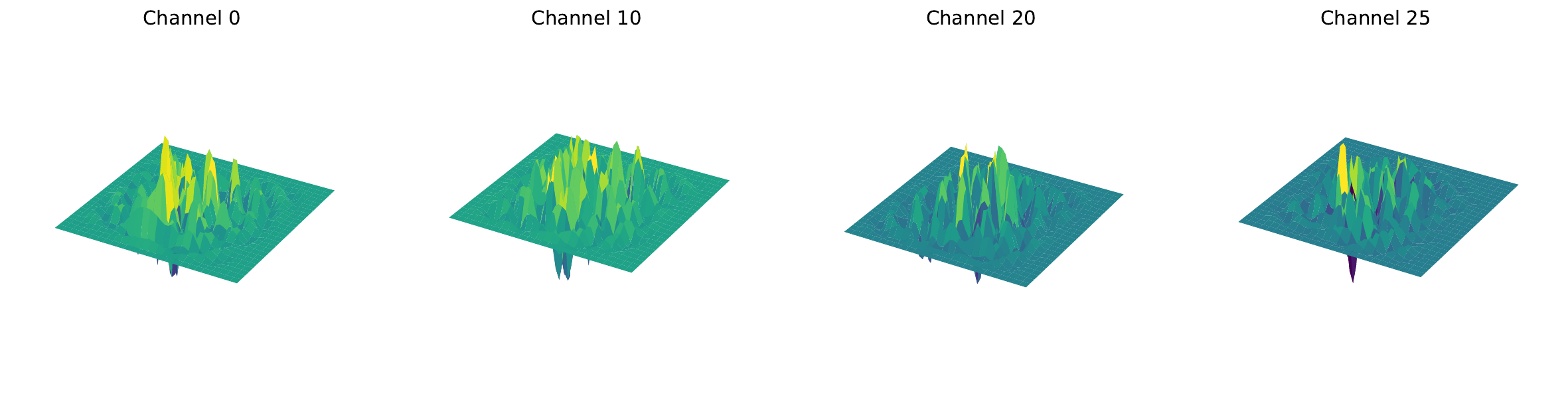}%
  }\hfill
  \subfloat[SVD spectrum across channels.\label{fig:channel_svd_spectrum}]{%
    \includegraphics[width=0.28\linewidth]{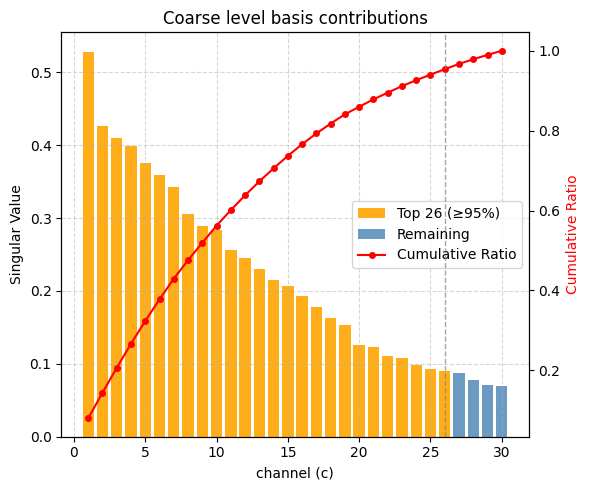}%
  }
  \caption{Diversity of the multi-channel coarse space. (a) Fine-grid responses obtained by activating channels $0$, $10$, $20$, and $25$ at the same coarsest-grid center and applying \eqref{eq:coarse_basis}; the responses are cropped around the heterogeneous scattering region. (b) Singular values of the matrix formed by vectorizing the per-channel responses; the gradual decay indicates that the channels span a non-redundant coarse subspace.}
  \label{fig:channel_diversity}
\end{figure}

This raises a natural question: are the channels genuinely doing independent work, or is the wide representation largely redundant? To answer it, we fix the center point of the coarsest grid, set one channel at this point to one and all other coarse-level entries to zero, and then map this tensor to the physical grid using \eqref{eq:coarse_basis}. Repeating this procedure over channels gives a family of fine-grid responses $\{\Psi^{(c)}\}_{c=0}^{C-1}$. Figure~\ref{fig:similarity_channel_basis_response} shows representative responses for channels $0$, $10$, $20$, and $25$, cropped around the heterogeneous scattering region. Although these responses are generated from the same spatial location and the same learned prolongation hierarchy, they exhibit different phase, orientation, and interference patterns. The channels therefore encode distinct physical components rather than replicated scalar bases.
We further quantify this diversity by stacking the vectorized responses $\Psi^{(c)}$ across all channels and computing their singular values. The resulting spectrum, shown in Figure~\ref{fig:channel_svd_spectrum}, decays gradually rather than collapsing to a few dominant modes. This indicates that the learned channels span a genuinely multi-dimensional coarse subspace. The multi-channel construction is therefore an essential part of the learned coarse representation, rather than a simple increase in feature width.

\subsection{Layer-by-Layer Progressive Finetuning for Large Domain Scalability}
\label{sec:llpf}
Achieving efficient convergence on large-scale computational domains presents a significant challenge for learning-based solvers. Models trained on small domains can accurately approximate local wave interactions but often fail to capture the long-range correlations essential for global error propagation. Consequently, when deployed on larger grids, these locally trained models exhibit degraded convergence rates, necessitating numerous outer iterations. A spectral interpretation of this long-range degradation is provided in Appendix~\ref{sec:app_spectral_analysis}.

While finetuning the entire model on the target large domain can mitigate this issue, it incurs prohibitive computational costs due to the increased grid resolution and the expense of backpropagating through the full depth of the multigrid hierarchy. To address this trade-off between scalability and efficiency, we propose the Layer-by-Layer Progressive Finetuning (LLPF) strategy. LLPF leverages the hierarchical nature of McMg to systematically extend the spatial support of the learned solver with minimal training overhead.

\begin{figure}[htbp]
    \centering
    \begin{tikzpicture}[
        node distance=0.15cm,
        block/.style={rectangle, draw=black, thick, minimum width=2.2cm, minimum height=0.7cm, align=center, font=\footnotesize\sffamily},
        active/.style={block, fill=cyan!20, draw=cyan!80!black},
        frozen/.style={block, fill=gray!20, text=gray!60, draw=gray!40, dashed},
        arrow/.style={->, >=Stealth, thick, color=black!70},
        label_t/.style={font=\bfseries\small, align=center}
    ]

    \node[active] (s1_l1) {Level 1 (Fine)};
    \node[active, below=of s1_l1] (s1_l2) {Level 2};
    \node[active, below=of s1_l2] (s1_l3) {Level 3};
    \node[active, below=of s1_l3] (s1_l4) {Level 4};
    \node[active, below=of s1_l4] (s1_l5) {Level 5 (Coarse)}; 
    
    \node[label_t, above=0.3cm of s1_l1] (title1) {Stage 1: Base Model\\ $\mathcal{M}_1$};
    \node[below=0.2cm of s1_l5, font=\scriptsize] (domain1) {Train on $\mathcal{D}_1$};

    \node[frozen, right=1.5cm of s1_l1] (s2_l1) {Level 1 (Fine)};
    \node[frozen, below=of s2_l1] (s2_l2) {Level 2};
    \node[frozen, below=of s2_l2] (s2_l3) {Level 3};
    \node[frozen, below=of s2_l3] (s2_l4) {Level 4};
    \node[frozen, below=of s2_l4] (s2_l5) {Level 5};
    \node[active, below=of s2_l5] (s2_l6) {Level 6 (New)};

    \node[label_t] at (title1 -| s2_l1) (title2) {Stage 2: Expansion\\ $\mathcal{M}_2$};
    \node[below=0.2cm of s2_l6, font=\scriptsize] (domain2) {Finetune on $\mathcal{D}_2$};

    \node[frozen, right=1.5cm of s2_l1] (s3_l1) {Level 1 (Fine)};
    \node[frozen, below=of s3_l1] (s3_l2) {Level 2};
    \node[frozen, below=of s3_l2] (s3_l3) {Level 3};
    \node[frozen, below=of s3_l3] (s3_l4) {Level 4};
    \node[frozen, below=of s3_l4] (s3_l5) {Level 5};
    \node[frozen, below=of s3_l5] (s3_l6) {Level 6};
    \node[active, below=of s3_l6] (s3_l7) {Level 7 (New)};

    \node[label_t] at (title1 -| s3_l1) (title3) {Stage 3: Expansion\\ $\mathcal{M}_3$};
    \node[below=0.2cm of s3_l7, font=\scriptsize] (domain3) {Finetune on $\mathcal{D}_3$};
    
    \path (s1_l2.east) -- (s2_l2.west) coordinate[midway] (mid1);
    \draw[arrow] ($(mid1)+(-0.4,0)$) -- ($(mid1)+(0.4,0)$);

    \path (s2_l3.east) -- (s3_l3.west) coordinate[midway] (mid2);
    \draw[arrow] ($(mid2)+(-0.4,0)$) -- ($(mid2)+(0.4,0)$);

    \node[anchor=north west, draw=black!20, rounded corners, fill=white, inner sep=5pt] (legend) 
            at ($(s1_l5.west |- s3_l7.south) + (0,-0.1)$) {
            \scriptsize
            \begin{tikzpicture}[baseline=(current bounding box.center)]
                \node[active, minimum width=0.5cm, minimum height=0.3cm] (la) at (0,0) {};
                \node[right=0.1cm of la] (txt_a) {Active}; 
                
                \node[frozen, minimum width=0.5cm, minimum height=0.3cm, right=0.6cm of txt_a] (lf) {};
                \node[right=0.1cm of lf] {Frozen};
            \end{tikzpicture}
        };

    \end{tikzpicture}
    \caption{Schematic of the Layer-by-Layer Progressive Finetuning (LLPF) strategy. The process begins with a base model $\mathcal{M}_1$ trained on a small domain $\mathcal{D}_1$. As the domain size increases, we append new coarser levels. During finetuning on larger domains ($\mathcal{D}_k$), parameters of the pre-existing fine levels (grey) are frozen to preserve local features, and only the newly added coarse operator (blue) is optimized to capture long-range interactions.}
    \label{fig:llpf_diagram}
\end{figure}
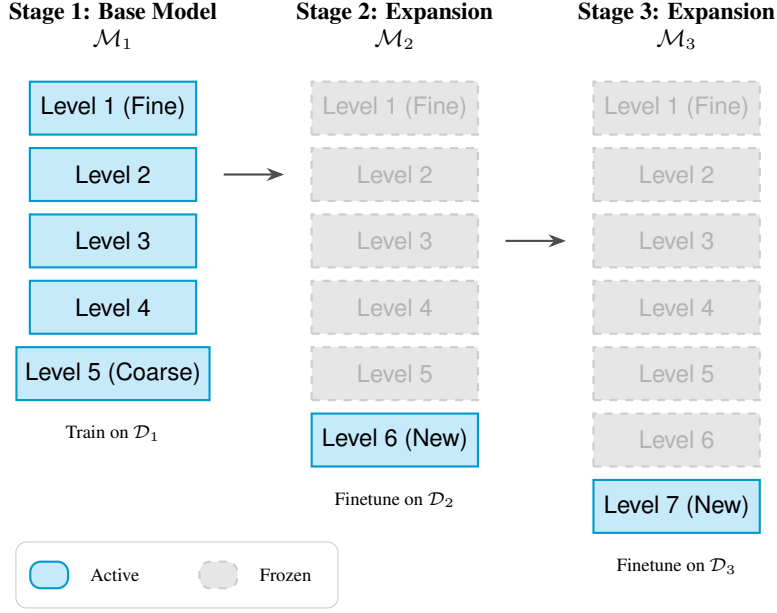
The LLPF procedure, illustrated in Figure~\ref{fig:llpf_diagram}, adapts a model trained on a small initial domain $\mathcal{D}_1$ to a sequence of increasingly larger domains $\{\mathcal{D}_k\}_{k=2}^K$, where $\mathcal{D}_K$ represents the target large-scale problem. The process consists of two main phases:

\begin{enumerate}
    \item \textbf{Base Model Initialization ($\mathcal{D}_1$):} 
    We commence by training a standard McMg solver with $L$ levels on the smallest domain $\mathcal{D}_1$. This step efficiently establishes the baseline parameters for resolving high-frequency local errors and short-range discrete physics. 
    
    \item \textbf{Recursive Expansion and Finetuning ($\mathcal{D}_k$):} 
    For each subsequent domain $\mathcal{D}_k$ larger than $\mathcal{D}_{k-1}$, we perform the following operations:
    \begin{itemize}
        \item \textbf{Architecture Growth:} We append an additional, coarser level to the bottom of the V-cycle hierarchy. This extends the effective receptive field of the solver, enabling it to represent the lower-frequency modes characteristic of the larger domain.
        \item \textbf{Progressive Locking:} The parameters of all previously trained levels (fine to intermediate grids) are frozen. These levels already capture the local physics correctly; preserving them prevents "catastrophic forgetting" and ensures stability.
        \item \textbf{Selective Optimization:} We finetune only the parameters associated with the newly added coarse level using data generated on $\mathcal{D}_k$. Since the fine-grid errors are already handled by the frozen layers, the optimization focuses exclusively on resolving the residual global error components.
    \end{itemize}
\end{enumerate}

This strategy keeps the training cost tractable as the problem size grows. By restricting backpropagation to a single variable level at each stage, LLPF extends the effective support of the learned Green's operator without requiring full-model retraining.

\section{Numerical Experiments}
In this section, we present a systematic evaluation of the proposed McMg method, benchmarking its performance against established baseline algorithms across diverse heterogeneous datasets. Unless otherwise stated, all computational frameworks were implemented using PyTorch 2.6.0 with CUDA 12.6 acceleration. The experiments were conducted on a high-performance workstation equipped with an Intel Xeon Gold 6444Y CPU and an NVIDIA RTX A6000 GPU.

\subsection{Benchmarking McMg against Classical and Neural Solvers} \label{sec:comp_with_other_model_OpenBreastUS}
We evaluate McMg across three challenging regimes, each stressing a different failure mode of classical solvers: a \emph{high-frequency} setting on the OpenBreastUS dataset~\cite{zeng2025openbreastus}, a \emph{high-contrast} setting on the CurveFault-B dataset, and a \emph{large-scale three-dimensional} setting on the Kimberlina 1.2 CCUS dataset, the latter two from OpenFWI~\cite{deng2022openfwi}. We deploy the learned McMg both as a standalone solver and as a preconditioner, and compare it against classical and neural solvers. We describe each setting and its training protocol below, then the evaluation protocol, and finally report training cost and solver performance.

\begin{figure}[htbp]
    \centering
    \subfloat[OpenBreastUS Breast dataset.\label{fig:OpenBreastUS_breast}]{%
        \includegraphics[width=0.48\linewidth]{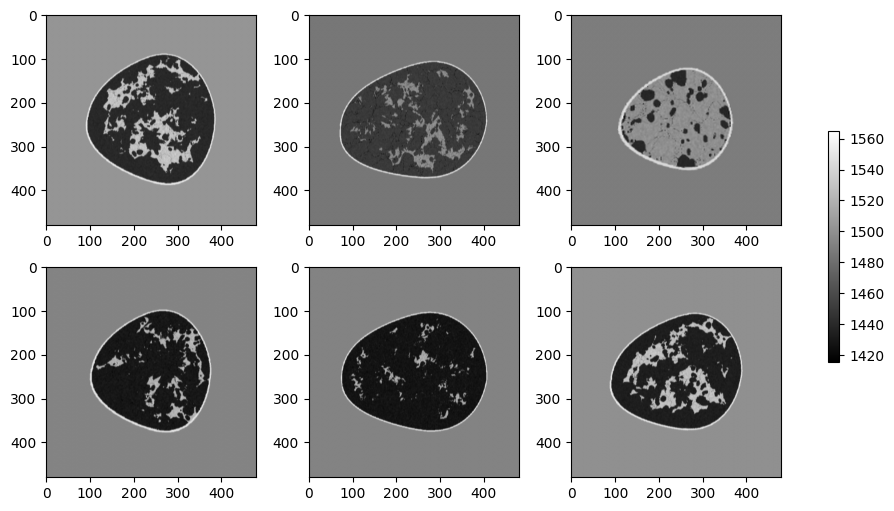}%
    }\hfill
    \subfloat[OpenFWI CurveFault-B dataset.\label{fig:openfwi_curvefault_B}]{%
        \includegraphics[width=0.48\linewidth]{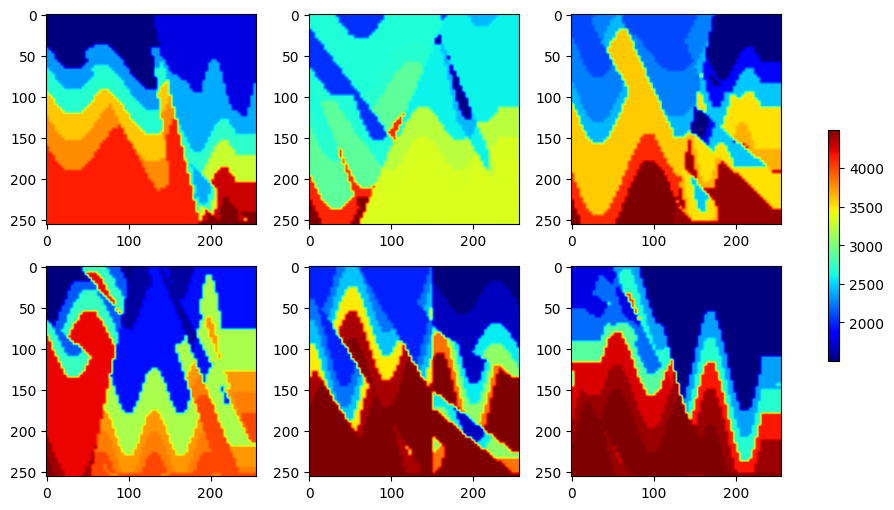}%
    }
    \caption{Representative velocity models for the two 2D benchmarks: (a) the high-frequency OpenBreastUS breast dataset and (b) the high-contrast OpenFWI CurveFault-B dataset.}
    \label{fig:datasets_2d}
\end{figure}

\paragraph{High-frequency regime (OpenBreastUS).}\label{sec:comp_with_CBS_OpenBreastUS}
The OpenBreastUS dataset comprises realistic 2D USCT breast models (Figure~\ref{fig:OpenBreastUS_breast}) on a $480\times 480$ grid. We target the challenging low points-per-wavelength ($\texttt{ppw}$) regime, where the pollution effect degrades the discrete operator's spectral properties and destabilizes iterative solvers; here McMg's multi-channel adaptive operators conform to local heterogeneity and remain stable where fixed-stencil approaches fail. Normalizing the domain to a unit square ($1\,\text{m}\times 1\,\text{m}$) gives a grid spacing $h = 1/480\,\text{m}$ and $\texttt{ppw} = 2\pi/(k_0 h) = 6$, where $k_0$ is the background wavenumber. 
The resulting dimensionless wavenumber $k_0 D \approx 500$ ($D$ the domain size) places the problem firmly in the high-frequency regime~\cite{livshits2006accuracy}. We use an absorbing layer of width 48 grid points and multigrid levels $[1,2,4,8,8]$, and train on 1000 samples for 400 epochs with batch size 10 using Adam.
\paragraph{High-contrast regime (CurveFault-B).}\label{sec:comp_with_other_openfwi}
The CurveFault-B dataset features sharp discontinuities and large velocity variations, from $c_{\min} = 1500$ m/s to $c_{\max} = 4500$ m/s (Figure~\ref{fig:openfwi_curvefault_B}). CBS is challenged here: higher contrasts enlarge the scattering potential and make the constant-background Green's preconditioner less contractive, leading to slower convergence. McMg instead generates adaptive stencils and smoothers with nonlinear dependence on the local wavenumber $k(x)$, capturing the wave physics at sharp interfaces and improving conditioning. We normalize the domain to the unit square $[0,1]\times[0,1]$ m on a $256\times 256$ grid, apply an absorbing layer of width 32 grid points, and set $\texttt{ppw} = 10$. We train on 3000 samples with batch size 30; all other settings follow the high-frequency configuration.

\paragraph{Large-scale 3D regime (Kimberlina 1.2 CCUS).}
Three-dimensional problems have long challenged AI-based PDE solvers, whose parameter counts and memory footprints grow rapidly with domain size. Thanks to its \emph{Setup} and \emph{Parameter Sharing} designs, McMg both trains and runs inference on $224\times 224\times 224$ volumes on a single NVIDIA A6000 GPU. We adopt the Kimberlina 1.2 CCUS dataset from OpenFWI, which models subsurface evolution across CO\textsubscript{2} injection years on an original $400\times 400\times 350$ grid. We interpolate each model to $192\times 192\times 192$ and add a 16-point absorbing layer, giving a $224\times 224\times 224$ domain with $\texttt{ppw}=6$, and sound speeds from $c_{\min}=2637$ to $c_{\max}=4201$ m/s---over 11 million degrees of freedom. We train on the first 35 years (630 samples) and validate on the 100th year (50 samples), using batch size 3 and multigrid levels $[1,2,4,6,8]$ with 20 channels; all other settings follow the high-frequency configuration. The resulting model has only 3.82 M parameters and uses roughly 45 GB of GPU memory for training and 6 GB for inference.

\paragraph{Baselines and evaluation protocol.}
On the high-frequency benchmark, we group our baselines into two classes. The classical solvers are CBS, a GMRES solver preconditioned by shifted-Laplacian multigrid (GMRES+SL)~\cite{erlangga2004class}, and a GPU-accelerated sparse LU factorization as a direct-solver reference; the neural solver is MGCFNN~\cite{xiemgcfnn}. We evaluate McMg as a standalone solver, as a GMRES preconditioner (GMRES+McMg), and as a Born-series preconditioner (Born+McMg). To ensure a fair comparison, all solvers act on the same five-point finite-difference discretization; the Born-series solvers (CBS and Born+McMg) retain their FFT-based Green operator but construct it from the discrete five-point symbol rather than the continuous Laplacian symbol, matching the discretization of the others (Appendix~\ref{sec:app_discrete_green}).
All iterative methods except CBS run in mixed precision, applying the preconditioner in FP32 and evaluating the residual in FP64; CBS runs entirely in FP64, as its preconditioning reduces to a single element-wise multiplication for which precision conversion would cost more than it saves. All runs share the same GPU environment, with iteration counts and wall-clock times averaged over 50 test samples at a relative tolerance of $1\times10^{-6}$. MGCFNN uses the same configuration as McMg. The sparse LU reference is implemented with the highly optimized NVIDIA cuDSS\footnote{\url{https://developer.nvidia.com/cudss}} library in single precision, and the convolutional forward pass is accelerated with PyTorch's just-in-time (JIT) compilation. The full parameterization of the classical baselines---the shifted-Laplacian shift, smoother, coarse solver, and GMRES restart length, together with the CBS background wavenumber and relaxation---is recorded in Appendix~\ref{sec:app_baseline_config}.

\paragraph{Training cost.}
Before turning to solver performance, we compare the training cost of McMg and MGCFNN on the high-frequency benchmark under an identical protocol (Table~\ref{tab:train_comp_mgcfnn_vs_mpmg}). McMg attains lower training and validation losses while using roughly an eighth of the parameters and a third of the per-epoch time, reflecting the efficiency of its parameter-sharing multigrid design.

\begin{table}[htbp]
\centering
\caption{Training comparison between MGCFNN and McMg under an identical protocol on the OpenBreastUS dataset.}
\label{tab:train_comp_mgcfnn_vs_mpmg}
\begin{tabular}{lcccc}
\toprule
Model & Params (MB) & Epoch Time (s) & Train Loss & Val Loss \\
\midrule
MGCFNN & 20.20 & 58.49 & $1.89 \times 10^{-2}$ & $2.57\times 10^{-2}$ \\
McMg   & \textbf{2.44} & \textbf{21.04} & $\mathbf{1.22 \times 10^{-2}}$ & $\mathbf{1.44\times 10^{-2}}$ \\
\bottomrule
\end{tabular}
\end{table}

\begin{figure}[htbp]
    \centering
    \includegraphics[width=1.\linewidth]{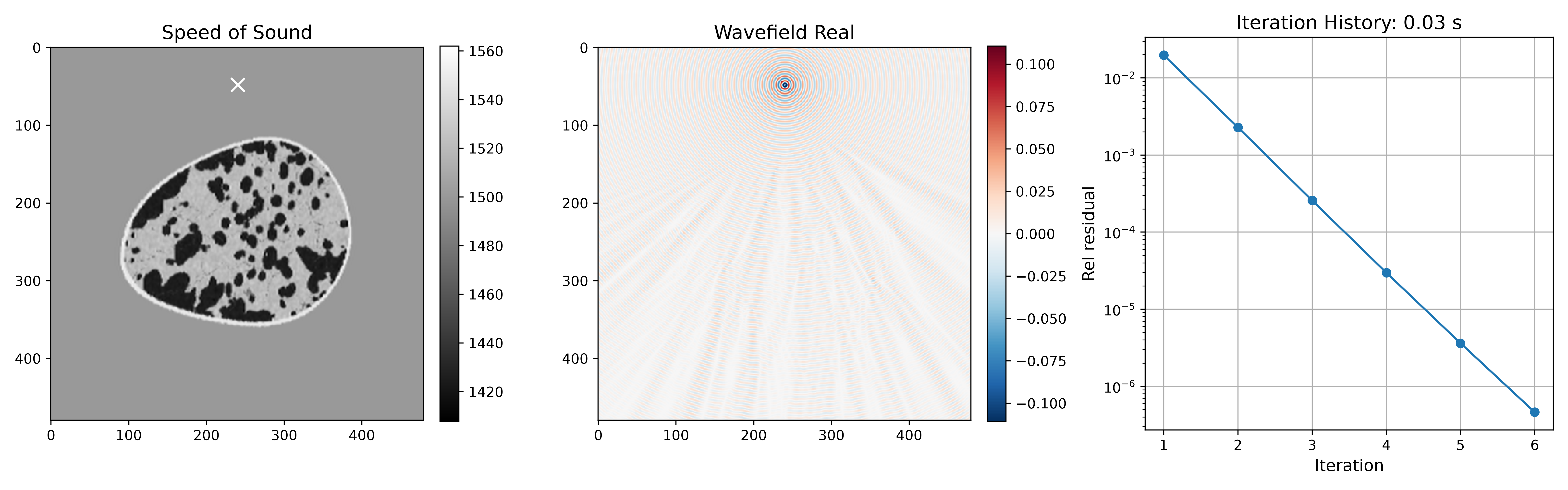}
    \caption{Inference performance of Born+McMg on the OpenBreastUS breast dataset, showing predicted wavefields and relative residual decay histories.}
    \label{fig:mpmg_bs_6005}
\end{figure}

\begin{figure}[htbp]
    \centering
    \includegraphics[width=1\linewidth]{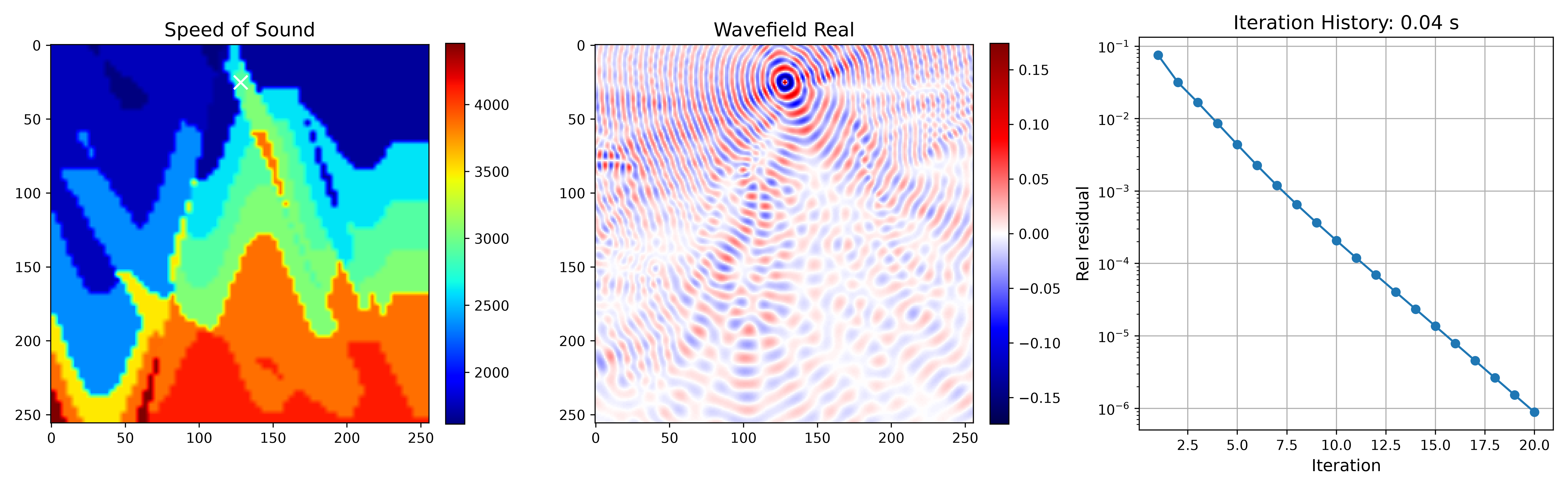}
    \caption{Inference performance of Born+McMg on the OpenFWI CurveFault-B dataset, showing predicted wavefields and relative residual decay histories.}
    \label{fig:openfwi_infer}
\end{figure}

\begin{figure}[htbp]
    \centering
    \includegraphics[width=1.\linewidth]{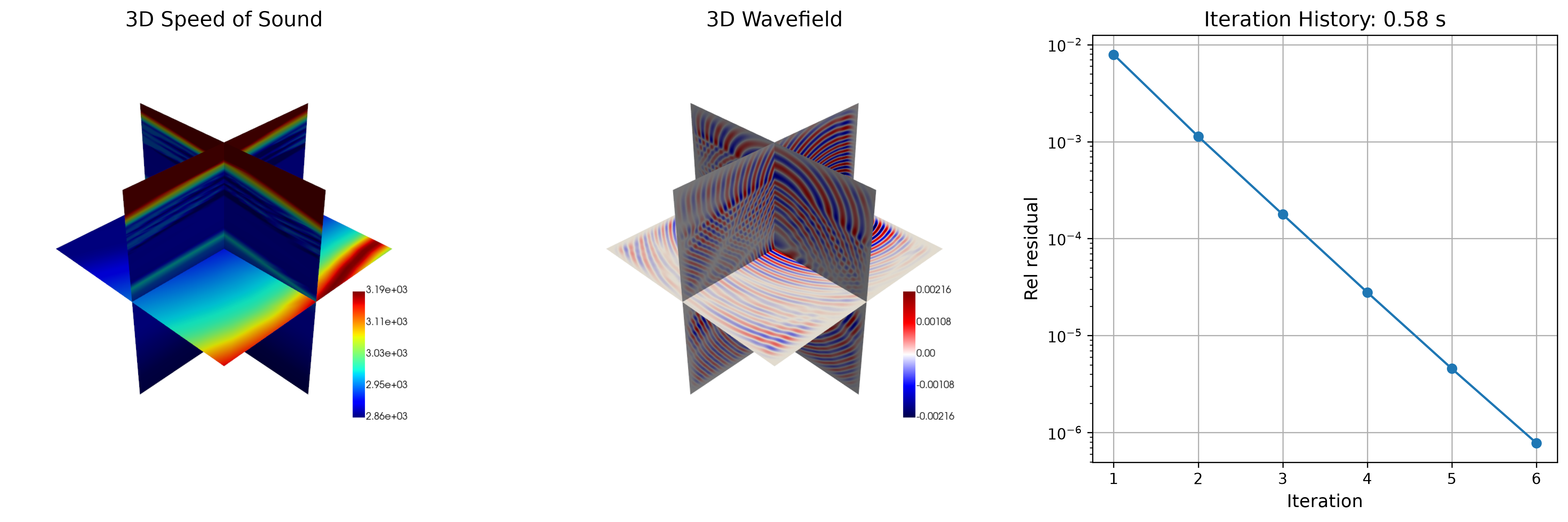}
    \caption{Inference performance of Born+McMg on the Kimberlina 1.2 CCUS dataset, showing predicted wavefields and relative residual decay histories.}
    \label{fig:3d_infer}
\end{figure}

\paragraph{Convergence and runtime.}
Table~\ref{tab:cbs_mpmg_precision} reports solver performance across the three regimes; we begin with the high-frequency benchmark. Among the classical methods, CBS attains the lowest wall-clock time ($0.372$ s), ahead of shifted-Laplacian-preconditioned GMRES (GMRES+SL, $3.290$ s) and a highly optimized sparse LU direct solver ($1.747$ s); we therefore adopt it as the reference classical solver throughout. Relative to CBS, McMg reduces the iteration count by $110.6\times$; because a single McMg iteration is more expensive than a single CBS iteration, this translates into a $10.3\times$ reduction in wall-clock time. That this gain originates in the learned preconditioner is confirmed by a controlled comparison: Born+McMg uses the same Born iteration as CBS and differs only in its preconditioner, yet converges in $8.0$ iterations as compared with $929.2$ for CBS. McMg likewise compares favorably with the neural baseline, converging nearly $4\times$ faster than MGCFNN. Interestingly, using McMg as a GMRES preconditioner reduces the iteration count further, from $8.4$ to $7.2$, but does not improve the wall-clock time ($0.036$ s versus $0.059$ s), indicating that the standalone McMg architecture is already efficient enough that the overhead of the Krylov recurrence outweighs the benefit of fewer iterations. The same advantage extends to the remaining regimes: relative to CBS, McMg reduces the iteration count by $51.3\times$ and $45.1\times$ and the wall-clock time by $3.4\times$ and $6.8\times$ in the high-contrast and three-dimensional settings, respectively, confirming its effectiveness across high-frequency, high-contrast, and large-scale regimes. Figures~\ref{fig:mpmg_bs_6005}, \ref{fig:openfwi_infer}, and \ref{fig:3d_infer} corroborate these findings, exhibiting accurate wavefields and rapid, monotone residual decay.
\begin{table}[htbp]
\centering
\caption{Solver performance across three benchmark regimes at a relative tolerance of $1\times10^{-6}$. Iteration counts and wall-clock times are averaged over 50 selected test samples; the best value in each group is shown in bold.}
\label{tab:cbs_mpmg_precision}
\renewcommand{\arraystretch}{1.15}
\begin{tabular}{@{}llcc@{}}
\toprule
Benchmark (Regime) & Method & Avg. Iters & Avg. Time (s) \\
\midrule
\multirow{7}{*}{\shortstack[l]{OpenBreastUS\\(high-frequency)}}
        & Sparse LU                     & --    & 1.747 \\
        & CBS                           & 929.2 & 0.372 \\
        & GMRES+SL                      & 781.5 & 3.290 \\
        & MGCFNN                        & 13.3  & 0.142 \\
        & McMg                          & 8.4   & \textbf{0.036} \\
        & GMRES+McMg                    & \textbf{7.2}   & 0.059 \\
        & Born+McMg                     & 8.0   & 0.037 \\
\midrule
\multirow{2}{*}{\shortstack[l]{CurveFault-B\\(high-contrast)}}
        & CBS                           & 934.2 & 0.131 \\
        & McMg                       & \textbf{18.2} & \textbf{0.038} \\
\midrule
\multirow{2}{*}{\shortstack[l]{Kimberlina 1.2 CCUS\\(3D large-scale)}}
        & CBS                           & 284.0 & 4.128 \\
        & McMg                          & \textbf{6.3} & \textbf{0.608} \\
\bottomrule
\end{tabular}

\end{table}


\subsection{Domain Generalization via Patch-Based Training}
The generalization capability of McMg is underpinned by the intrinsic locality of the learned stencils. Since the stencil coefficients at a given point are primarily governed by the local wave velocity structure and discretization errors rather than by the total number of grid points, the learned map can transfer from small patches to larger domains. This property enables a train-on-patches, deploy-on-full-domains strategy, where McMg is trained on computationally inexpensive subdomains and then tested on larger domains without finetuning.

\begin{figure}[htbp]
    \centering
    \begin{minipage}[b]{0.48\linewidth}
        \centering
        \includegraphics[width=\linewidth, keepaspectratio]{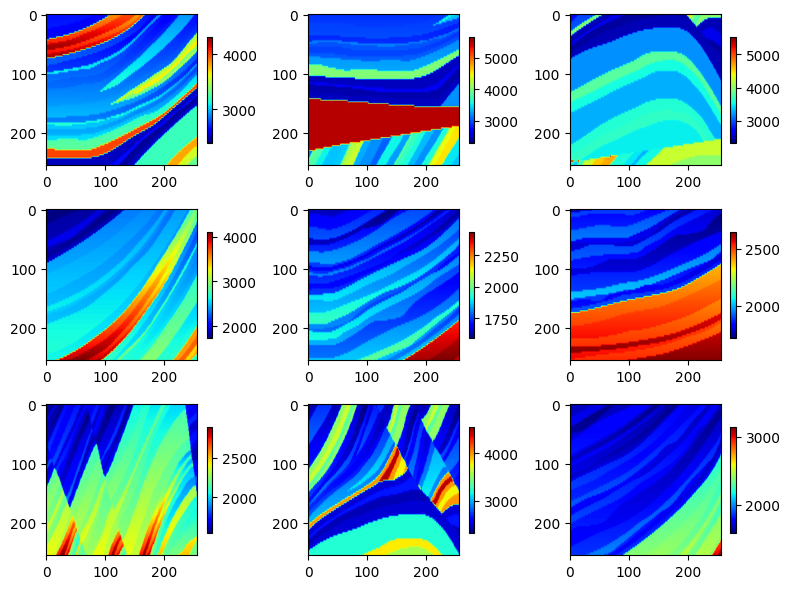}
        \\[0.5em] 
        \centerline{(a)}
    \end{minipage}
    \hfill
    \begin{minipage}[b]{0.48\linewidth}
        \centering
        \includegraphics[width=\linewidth, keepaspectratio]{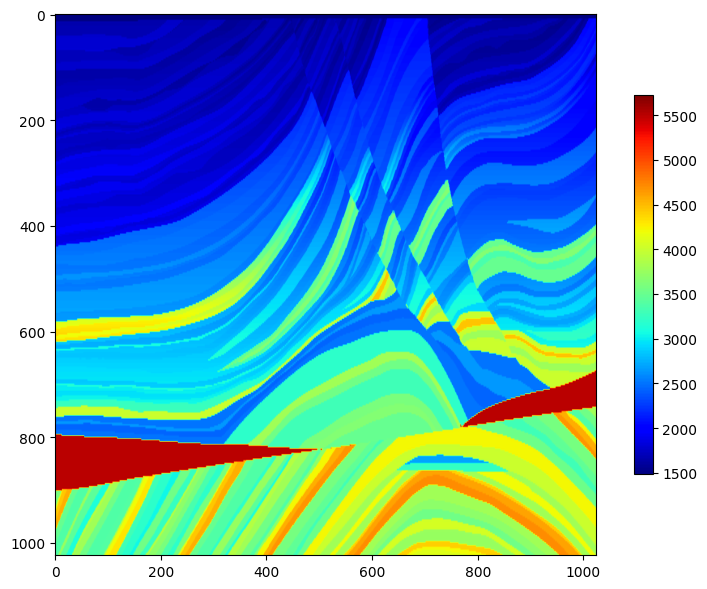}
        \\[0.5em]
        \centerline{(b)}
    \end{minipage}
    \caption{Velocity models used for generalization analysis: (a) a representative $256 \times 256$ patch randomly cropped for training, and (b) the full $1024 \times 1024$ Marmousi model used for testing.}
    \label{fig:marmousi_data}
\end{figure}

To validate this capability, we train McMg exclusively on small sub-patches ($256 \times 256$) randomly extracted from the Marmousi velocity model (Figure~\ref{fig:marmousi_data}a). Subsequently, the trained model is employed to solve the Helmholtz equation on the complete $1024 \times 1024$ domain (Figure~\ref{fig:marmousi_data}b). The experimental setup follows the configuration in Section~\ref{sec:comp_with_other_openfwi}, utilizing a sampling density of $\texttt{ppw} = 12$. The equation is discretized using a standard 5-point FDM, augmented with a sponge layer of width 32 to suppress boundary reflections.

\begin{figure}[htbp!]
    \centering
    \includegraphics[width=0.9\linewidth]{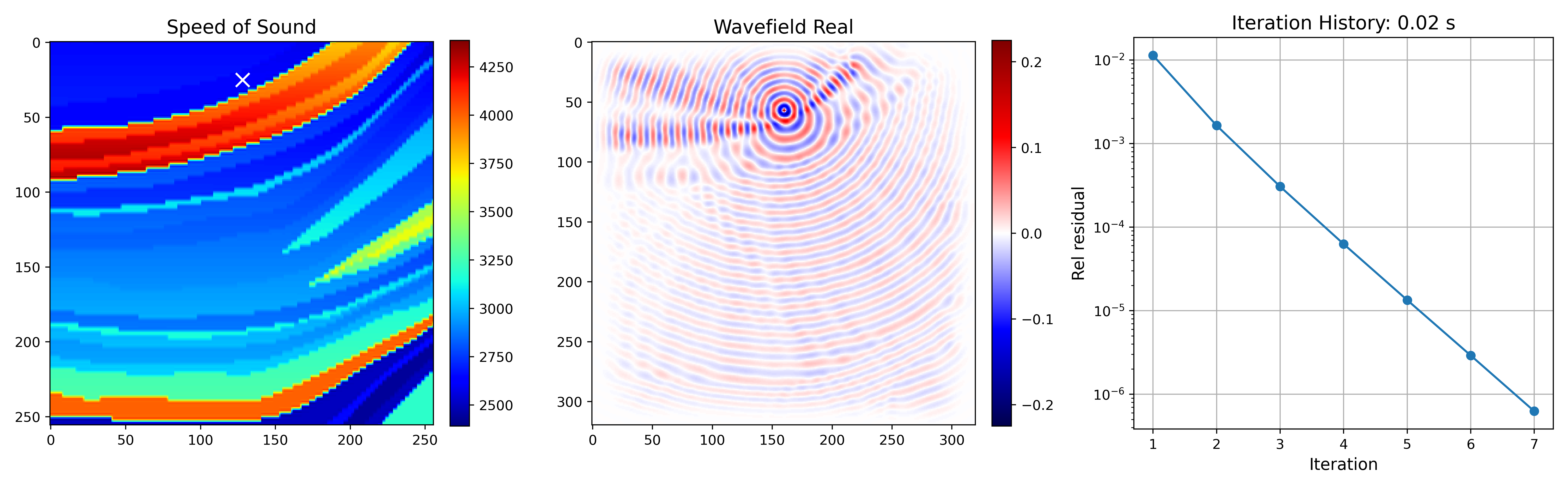}
    \caption{Convergence performance of McMg on a randomly cropped Marmousi patch.}
    \label{fig:mpmg_marmousi_part}
\end{figure}

\begin{figure}[htbp!]
    \centering
    \includegraphics[width=0.9\linewidth]{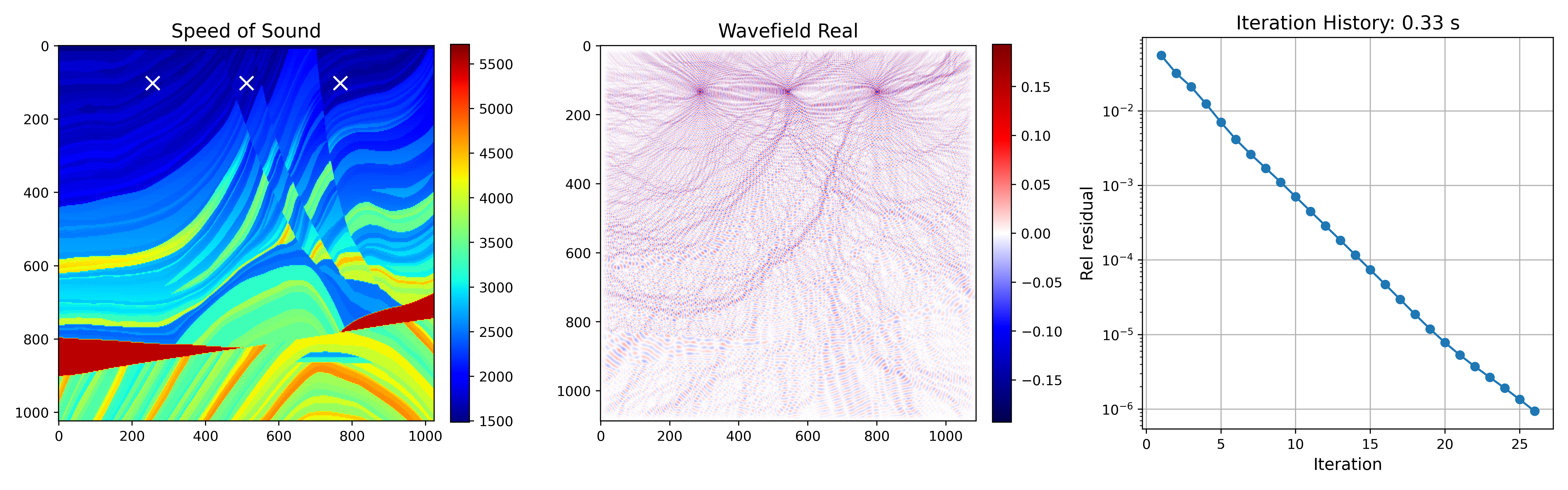}
    \caption{Convergence performance of McMg on the full Marmousi model.}
    \label{fig:mpmg_marmousi}
\end{figure}

Figure~\ref{fig:mpmg_marmousi_part} and Figure~\ref{fig:mpmg_marmousi} illustrate the solver's convergence performance on a cropped validation patch and the full domain, respectively. As anticipated, the model demonstrates robust generalization: the number of iterations required for the full-scale problem exhibits only a marginal increase compared to the sub-problems. This confirms that the local operators learned by McMg effectively capture the underlying physics, regardless of the global scale of the simulation.

\subsection{Scalability and Cross-Scale Generalization}\label{sec:scalability_large_domain}
A pivotal advantage of our McMg is the robust generalization across scales and wavenumbers. Our model, trained on small-scale domains, can be directly applied to large-scale problems with significantly higher wavenumbers without retraining. This ``train on small-scale, deploy on large-scale'' capability represents a substantial advancement over neural operators that typically require retraining when the problem size or frequency range changes. We investigate this scalability by training McMg on a limited-scale dataset ($480\times480$) and subsequently deploying it on problems of exponentially increasing size and frequency, ranging from a grid size of $512 \times 512$ (with $k_0 D \approx 500$) up to $8192 \times 8192$ (with $k_0 D \approx 8000$). Based on these results, we further improve large-domain scalability using the LLPF strategy described in Section~\ref{sec:llpf}. The experiments are conducted on the OpenBreastUS Breast dataset, utilizing the same configuration detailed in Section~\ref{sec:comp_with_other_model_OpenBreastUS} with a convergence tolerance of $10^{-6}$.

To evaluate scalability, we compare McMg, Born+McMg, and the CBS solver (Table~\ref{tab:scalability_comp}). For McMg, the per-iteration cost scales linearly as $\mathcal{O}(N)$ with respect to the total number of grid points $N$. For Born+McMg and CBS, which apply the Helmholtz operator via the FFT, the per-iteration cost is $\mathcal{O}(N\log N)$: the $\log N$ factor stems solely from the FFT-based operator application, while all remaining operations scale linearly in $N$.

\begin{table}[htbp]
\centering
\caption{Scalability comparison of McMg, Born+McMg, and CBS across increasing wavenumbers and grid sizes, with a single training on a small domain ($480\times480$).}
\label{tab:scalability_comp}

\resizebox{0.8\textwidth}{!}{%
\begin{tabular}{cc ccc ccc}
\toprule
\multirow{2}{*}{$k_0 D$} & \multirow{2}{*}{Grid size} & \multicolumn{3}{c}{Avg. Iters} & \multicolumn{3}{c}{Avg. Time (s)} \\
\cmidrule(lr){3-5} \cmidrule(lr){6-8}
 & & McMg & Born+McMg & CBS & McMg & Born+McMg & CBS \\
\midrule
500  & 512$\times$512   & 8.7 & 6.3  & 860.2  & 0.04 & 0.03 & 0.47 \\
1000 & 1024$\times$1024 & 8.8 & 7.0  & 1011.0 & 0.12 & 0.09 & 1.74 \\
2000 & 2048$\times$2048 & 12.4 & 11.0 & 1825.8 & 0.61 & 0.53 & 11.5 \\
4000 & 4096$\times$4096 & 22.0 & 22.0 & 3603.4 & 4.07 & 5.59 & 118.0 \\
\bottomrule
\end{tabular}%
}
\end{table}

\begin{figure}[htbp]
  \centering
  \subfloat[McMg pretrained on a $512\times512$ domain.\label{fig:green_scalability:a}]{%
    \includegraphics[width=0.95\linewidth]{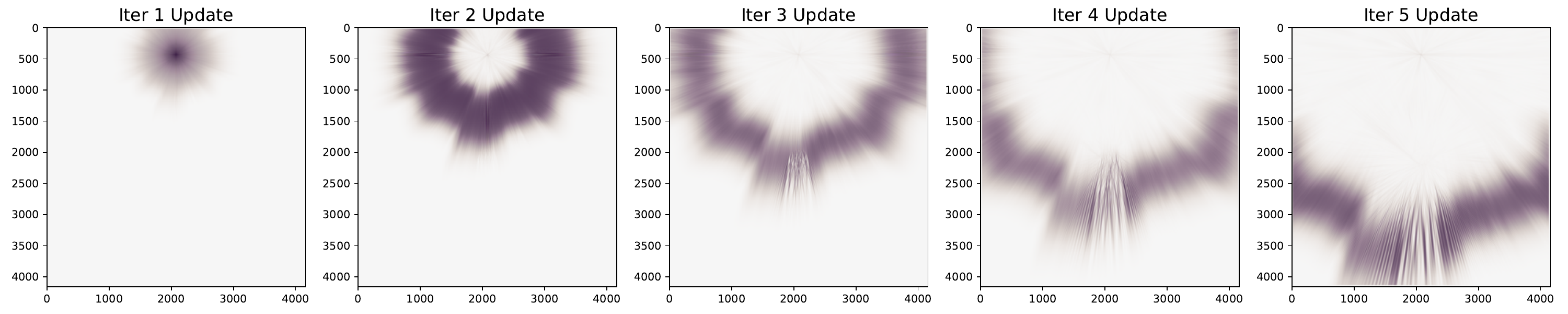}%
  }\\[0.8em]
  
  \subfloat[McMg finetuned on a $1024\times1024$ domain.\label{fig:green_scalability:b}]{%
    \includegraphics[width=0.95\linewidth]{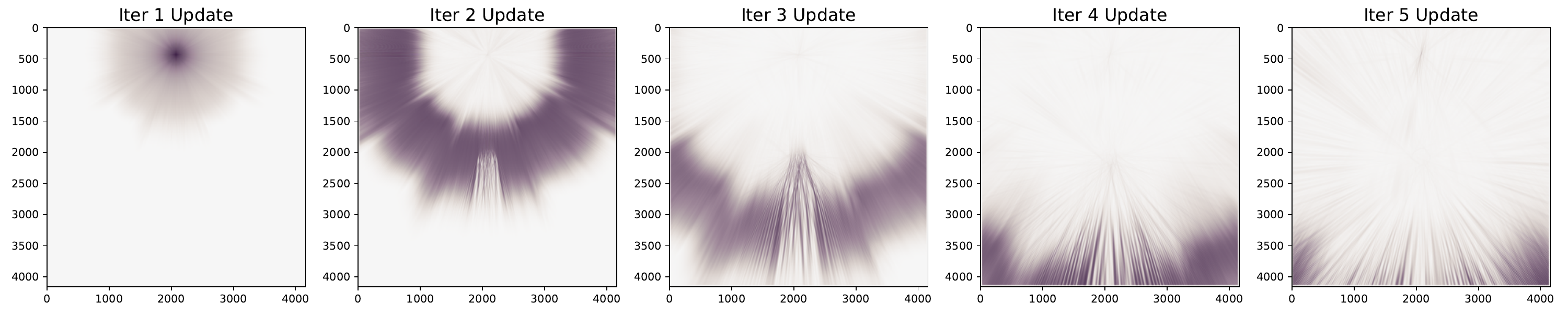}%
  }\\[0.8em]
  
  \subfloat[McMg finetuned with LLPF on a $1024\times1024$ domain.\label{fig:green_scalability:c}]{%
    \includegraphics[width=0.95\linewidth]{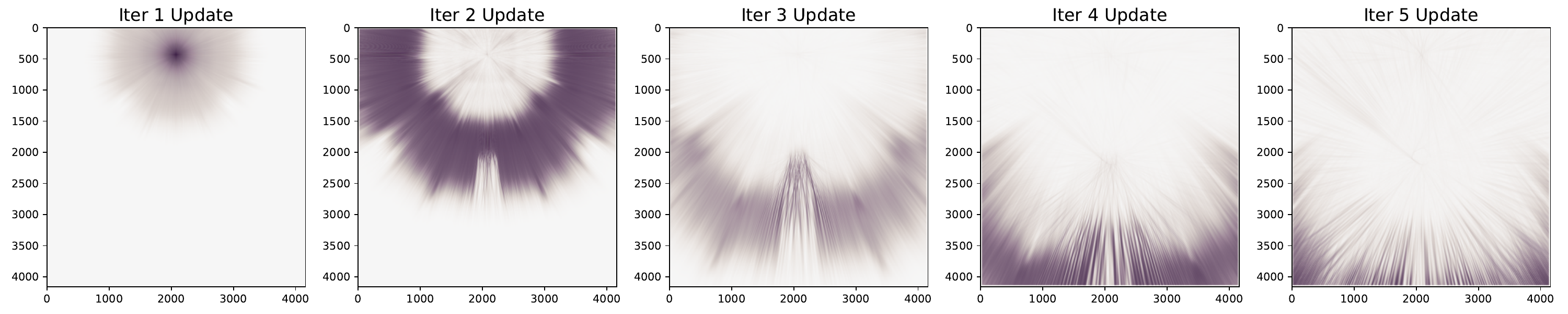}%
  }

  \caption{Effective support of the learned Green's function on a $4096\times4096$ domain, probed with a point source. Each row shows the first five corrections $\Delta\bm u^{(k)}=\mathcal{MG}(\bm r_{\mathrm{H}}^{(k)})$, $k=0,\dots,4$; by \eqref{eq:green_visual} the first panel approximates the discrete Green's function $A^{-1}\bm\delta_{x_s}$. (a) Pretraining on $512\times512$ yields a spatially truncated, short-range correction. (b) Direct finetuning on $1024\times1024$ and (c) LLPF on $1024\times1024$ both extend the correction over a much larger range, with LLPF matching (b) at a fraction of the finetuning cost.}
  \label{fig:green_scalability}
\end{figure}

The iteration count nonetheless increases monotonically with the domain size. To see why, recall the McMg correction step \eqref{eq:mcmg_iteration}, so the correction applied at step $k$ is $\Delta\bm u^{(k)} = \mathcal{MG}(\bm r_{\mathrm{H}}^{(k)})$. Driving the solver with a point source $\bm f = \bm\delta_{x_s}$ from a zero initial guess, the first correction
\begin{equation}\label{eq:green_visual}
    \Delta\bm u^{(0)} = \mathcal{MG}(\bm\delta_{x_s}) \approx A^{-1}\bm\delta_{x_s}
\end{equation}
is exactly the learned approximation of the discrete Green's function centered at the source $x_s$, and each subsequent correction applies the same operator to the remaining residual, transporting the wave further across the domain. Visualizing the first five corrections $\{\Delta\bm u^{(k)}\}_{k=0}^{4}$ (Figure~\ref{fig:green_scalability}) therefore exposes the effective spatial support of the learned Green's function. A model trained only on a small domain reproduces a spatially truncated Green's function (Figure~\ref{fig:green_scalability:a}): it acts as a localized operator that captures only short-range interactions, so each correction reaches a limited range and many outer iterations are needed to propagate the solution across a larger domain. To mitigate this growth in iteration count, we next consider finetuning and the LLPF strategy.

We evaluate different finetuning strategies, with results reported in Table~\ref{tab:mpmg_bs_scalability_comp_finetune_and_llpf}. Direct finetuning on the $1024 \times 1024$ domain with $7$ levels yields the best overall performance, balancing optimization stability and effective Green's function support. In contrast, finetuning directly on the $2048 \times 2048$ domain often degrades convergence due to increased optimization difficulty (see Appendix~\ref{sec:app_spectral_analysis}). LLPF mitigates this issue by freezing existing levels and finetuning only the newly added coarsest level. As shown in Table~\ref{tab:mpmg_bs_scalability_comp_finetune_and_llpf}, LLPF consistently improves convergence as coarse levels are added, matching the convergence of direct finetuning on the $1024 \times 1024$ domain at a small fraction of its cost. Because backpropagation is confined to the single newly added level, extending a five-level model to six levels with LLPF takes only $0.88$h, about $3\times$ faster than directly finetuning a full six-level model on the same $1024\times1024$ domain ($2.51$h), with a comparable reduction in memory. The effective support of the learned Green's function tells the same story (Figure~\ref{fig:green_scalability}): although LLPF updates only the added level, its correction (Figure~\ref{fig:green_scalability:c}) attains the same long-range support as direct finetuning of the full model on the $1024\times1024$ domain (Figure~\ref{fig:green_scalability:b}), and both extend markedly beyond the short-range, truncated correction of the model pretrained on $512\times512$ (Figure~\ref{fig:green_scalability:a}).

\begin{table}[htbp]
\centering
\caption{Comparison of pretraining (P.T.), direct finetuning (F.T), and LLPF for McMg across increasing domain sizes. Iteration counts are reported for different grid sizes and multigrid levels. LLPF progressively adds coarse levels while freezing previously trained parameters. Columns are color-coded by the training domain size: \colorbox{regfivetwelve}{$512$}, \colorbox{regtenfour}{$1024$}, and \colorbox{regtwentyfour}{$2048$}. Training/Finetuning time (P.T./F.T. Time) is measured on an NVIDIA H100 GPU over 400/200 epochs.}
\label{tab:mpmg_bs_scalability_comp_finetune_and_llpf}

\resizebox{0.98\textwidth}{!}{%
\begin{tabular}{cc
  >{\columncolor{regfivetwelve}}c >{\columncolor{regfivetwelve}}c >{\columncolor{regfivetwelve}}c
  >{\columncolor{regtenfour}}c >{\columncolor{regtenfour}}c >{\columncolor{regtenfour}}c
  >{\columncolor{regtwentyfour}}c >{\columncolor{regtwentyfour}}c >{\columncolor{regtwentyfour}}c
  >{\columncolor{regtenfour}}c >{\columncolor{regtwentyfour}}c}
\toprule
\multirow{2}{*}{$k_0 D$} & \multirow{2}{*}{Grid size}
& \multicolumn{3}{c}{P.T. 512}
& \multicolumn{3}{c}{F.T. 1024}
& \multicolumn{3}{c}{F.T. 2048}
& \multicolumn{2}{c}{LLPF} \\
\cmidrule(lr){3-5} \cmidrule(lr){6-8} \cmidrule(lr){9-11} \cmidrule(lr){12-13}
& & 5 lvls & 6 lvls & 7 lvls
& 5 lvls & 6 lvls & 7 lvls
& 5 lvls & 6 lvls & 7 lvls
& 6 lvls & 7 lvls \\
\midrule
500  & $512\times512$     & 7.0  & 6.7  & 6.8  & 10.4 & 7.5  & 7.7  & 18.5 & 11.6 & 11.5 & 7.0  & 7.1  \\
1000 & $1024\times1024$   & 7.2  & 7.0  & 7.0  & 11.0 & 7.3  & 8.3  & 27.3 & 16.0 & 15.2 & 7.2  & 7.2  \\
2000 & $2048\times2048$   & 11.0 & 10.9 & 10.0 & 14.0 & 8.0  & 8.3  & 42.5 & 18.9 & 19.9 & 8.3  & 8.5  \\
4000 & $4096\times4096$   & 22.0 & 21.1 & 19.0 & 23.9 & 12.0 & 11.1 & 69.0 & 25.0 & 29.0 & 13.0 & 12.3 \\
8000 & $8192\times8192$   & 53.0 & 55.0 & 43.0 & 42.0 & 21.0 & 20.0 &116.0 & 42.0 & 40.0 & 23.0 & 22.0 \\
\midrule
\multicolumn{2}{c}{\textbf{P.T./F.T. Time (h)}}
& $1.23$ & $1.28$ & $1.33$
& $2.43$ & $2.51$ & $2.55$
& $8.61$ & $8.76$ & $8.85$
& $\mathbf{0.88}$ & $\mathbf{3.0}$ \\
\bottomrule
\end{tabular}%
}
\end{table}

\subsection{Spectral Convergence Analysis}
We analyze the error evolution in the frequency domain to investigate how McMg handles spectral error modes. A major limitation of GMG is its inability to attenuate characteristic components~\cite{cui2025neural} (error modes with wavenumber $|\mathbf{\xi}| \approx k$). These modes are difficult to eliminate because they generate small residuals for local smoothers and suffer from severe phase errors on scalar coarse grids.

\begin{figure}[htbp]
    \centering
    \subfloat[]{%
        \includegraphics[width=0.32\linewidth]{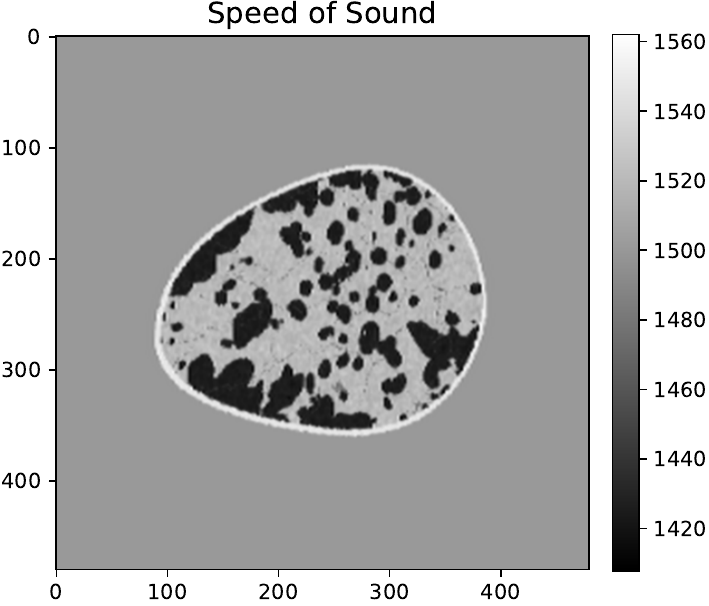}%
        \label{fig:1a}%
    }\hfill
    \subfloat[]{%
        \includegraphics[width=0.32\linewidth]{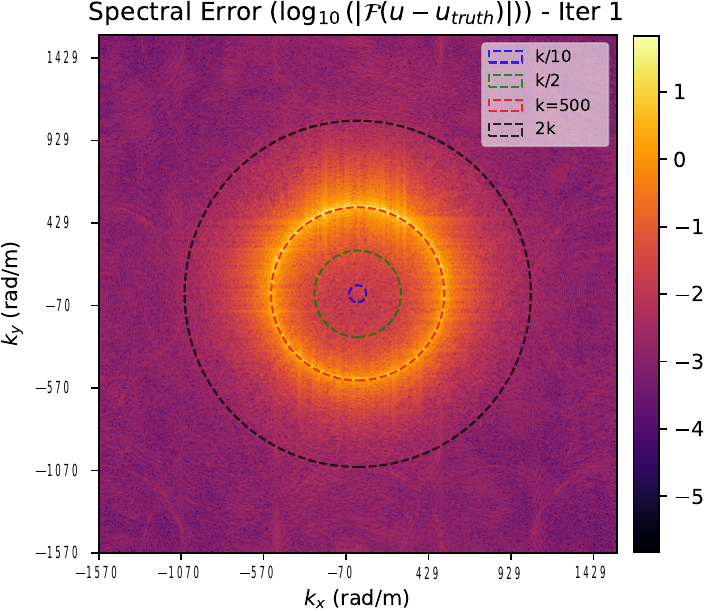}%
        \label{fig:1b}%
    }\hfill
    \subfloat[]{%
        \includegraphics[width=0.32\linewidth]{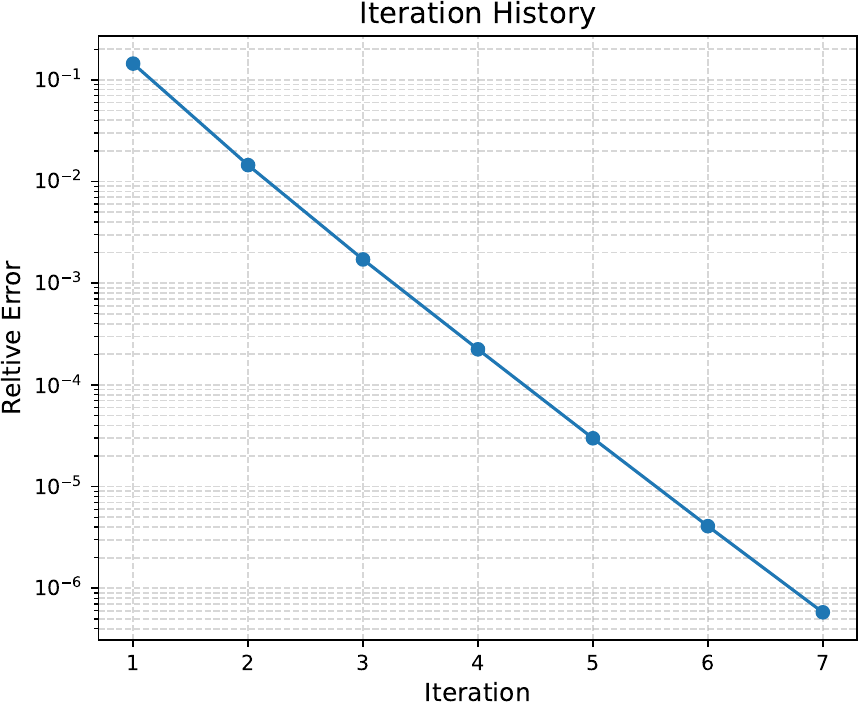}%
        \label{fig:1c}%
    }

    \vspace{0.5em}

    \subfloat[]{%
        \includegraphics[width=0.32\linewidth]{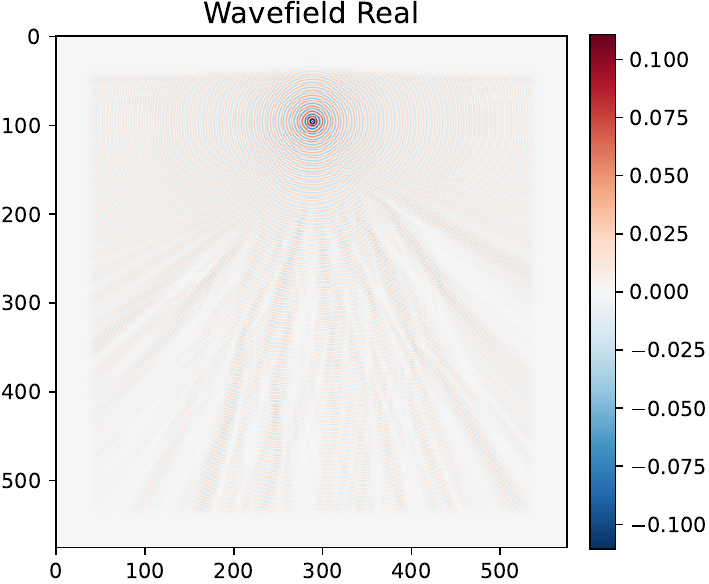}%
        \label{fig:1d}%
    }\hfill
    \subfloat[]{%
        \includegraphics[width=0.32\linewidth]{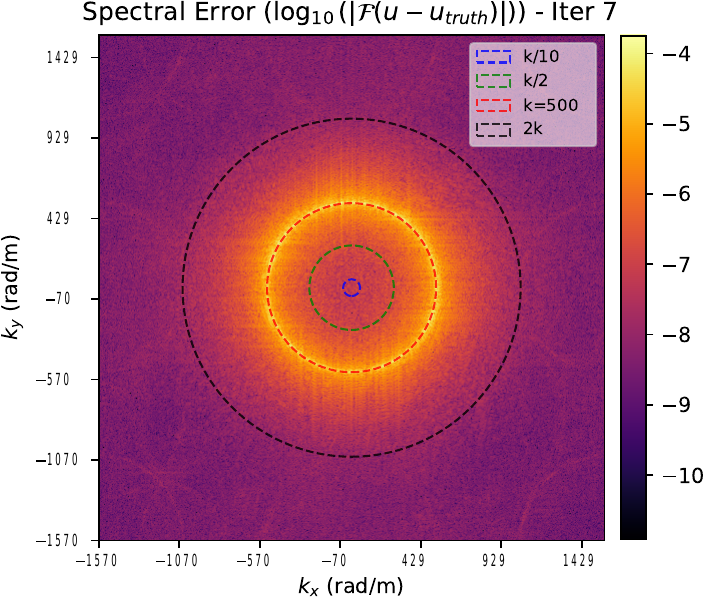}%
        \label{fig:1e}%
    }\hfill
    \subfloat[]{%
        \includegraphics[width=0.32\linewidth]{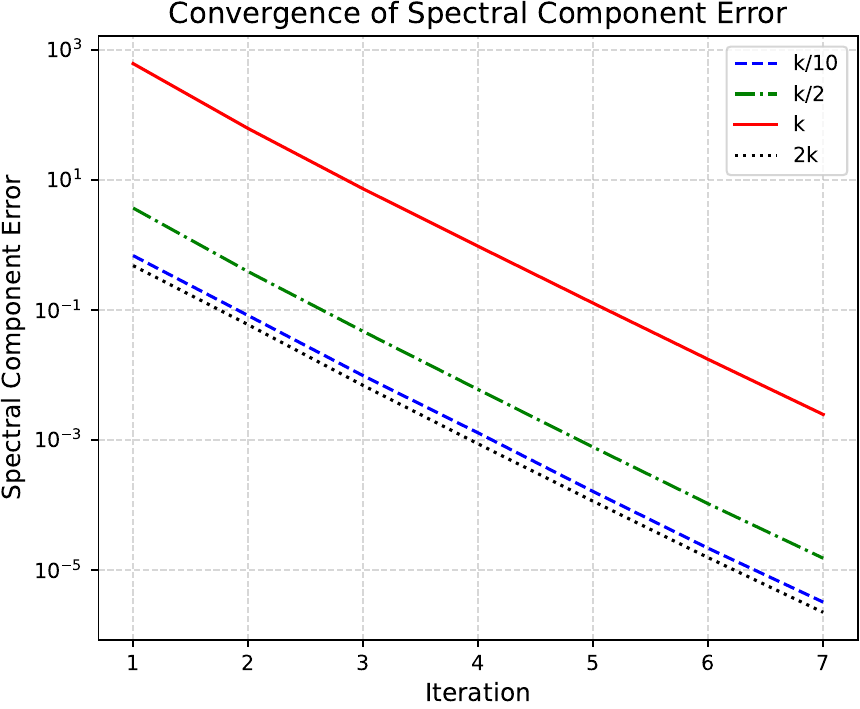}%
        \label{fig:1f}%
    }

    \caption{Spectral convergence analysis. (a) Sound speed map. (d) Wavefield real part. (b, e) Spectral error at iterations 1 and 7; dashed rings mark wavenumbers $k/10$, $k/2$, $k$, and $2k$. Note the initial error concentration near the characteristic mode $|\mathbf{\xi}| \approx k$ in (b). (f) Component-wise convergence. The difficult characteristic mode ($k$, red line) decays at a rate comparable to the lower- and higher-frequency bands.}
    \label{fig:all}
\end{figure}

We examine a test sample with a heterogeneous breast sound speed map (Figure~\ref{fig:1a}). Figures~\ref{fig:1b} and \ref{fig:1e} visualize the spectral error magnitude at iterations 1 and 7, respectively. The initial error (Figure~\ref{fig:1b}) is heavily concentrated around the characteristic ring $|\mathbf{\xi}| \approx k$, confirming that these modes dominate the error landscape. McMg suppresses these components by iteration 7 (Figure~\ref{fig:1e}), while also reducing the lower- and higher-frequency bands.

The quantitative breakdown in Figure~\ref{fig:1f} confirms this observation. The error at the characteristic frequency $k$ decays at an approximately exponential rate, parallel to both the low-frequency ($k/10$) and high-frequency ($2k$) components. This comparable decay across diverse bands supports the interpretation that McMg's multi-channel coarse representation and adaptive operators preserve the wave information needed to reduce characteristic modes that typically stagnate standard GMG solvers.

\section{Conclusions}
We have presented McMg, a learned multigrid preconditioner for heterogeneous Helmholtz equations that replaces end-to-end solution prediction with iterative residual-to-correction learning. The main architectural idea is to replace scalar coarse-grid correction by a multi-channel phase-space coarse representation: physical space is coarsened, but each coarse node carries latent wave-packet coefficients that can represent local phase, direction, and scattering information. Medium-dependent neural PDE operators and smoothers then adapt this packet representation to the local wave speed.
Two design choices make the method practical at scale. First, parameter sharing within the multigrid hierarchy yields a compact model. Second, the setup phase precomputes and caches the nonlinear physical features once per medium, so that each iteration reduces to a sequence of linear convolutions with cached coefficients. McMg can also be incorporated into the CBS framework, yielding Born+McMg. This variant is particularly effective in strongly scattering media, where training in the preconditioned Born-series residual space improves stability and convergence efficiency. To extend the solver to large domains without full retraining, LLPF preserves previously trained levels and finetunes only newly appended coarse levels, improving large-domain convergence at much lower training cost.

Across high-frequency (OpenBreastUS), high-contrast (OpenFWI CurveFault-B), and large-scale three-dimensional (Kimberlina) benchmarks, McMg converges in substantially fewer iterations and shorter wall-clock time than the fastest traditional solver in our tests, and outperforms recent neural solvers at comparable accuracy. Relative to CBS, this strongest traditional reference, McMg reduces iteration counts by $110.6\times$, $51.3\times$, and $45.1\times$ on the three benchmarks, with corresponding wall-clock speedups of $10.3\times$, $3.4\times$, and $6.8\times$. Spectral and coarse-space analyses further indicate that the learned operators reduce characteristic error modes near $|\bm\xi|\approx k$ and, by attaching a multi-channel local basis to each coarse node, represent oscillatory components on grids where a scalar coarse space would be too sparse.

Several directions remain open. A local Fourier or Wigner diagnostic should directly test whether the learned channels tile the characteristic shell $|\xi|\approx k(x)$; effective coarse-space approximation tests should compare McMg against scalar GMG, LOD, and GMsFEM bases; and stability diagnostics should report field-of-values, pseudospectral, or worst-case residual amplification measures for the learned preconditioner. The present formulation relies on structured grids, so extending the learned stencils to unstructured meshes and complex geometries is important. Finally, McMg should be embedded as an inner solver in large-scale inverse problems, where repeated Helmholtz solves for many sources and frequencies are a dominant computational cost.

\bibliographystyle{plainnat}
\bibliography{ref_arxiv}
\appendix
\section{Additional Derivations and Experiments}

\subsection{Alternative Discretizations and Boundary Conditions}
\label{sec:app_discretization}
The main text discretizes the Helmholtz operator with the five-point finite-difference (FDM) stencil~\eqref{eq:FDM_stencil} and truncates the domain with a sponge layer. This appendix details the two principal alternatives, the perfectly matched layer (PML) as an absorbing boundary and the Fourier spectral method (FSM) as a discretization, and then compares all four combinations empirically.

\subsubsection{Derivation of the Perfectly Matched Layer (PML)}
\label{sec:pml_derivation}
The Perfectly Matched Layer (PML) method extends the spatial coordinates into the complex plane to introduce artificial decay without reflection. For the Helmholtz equation, this is achieved by applying a complex coordinate stretching:  
\begin{equation}
    \tilde{x}_j = x_j + i p_j(x_j), \quad j = 1, \dots, D,
\end{equation}
where $p_j(x_j)$ is a monotonically increasing function in the $j$-th coordinate direction. The corresponding Jacobian of the transformation is defined as
\begin{equation}
    \frac{\partial \tilde{x}_j}{\partial x_j} = 1 + i \frac{\partial p_j(x_j)}{\partial x_j}, \quad j = 1, \dots, D.
\end{equation}
Introducing the complex stretching factor:
\begin{equation}
    s_j(x_j) = 1 + i \frac{\partial p_j(x_j)}{\partial x_j} = 1 + i\sigma_j(x_j),
\end{equation}
where $\sigma_j(x_j)$ denotes the damping profile in the $j$-th direction, we have $\sigma_j(x_j) = 0$ inside the physical domain $\Omega$ and $\sigma_j(x_j) > 0$ within the absorbing layer $\Omega_{\text{abl}}$ (often chosen as a smooth polynomial). The derivative with respect to the complex coordinate becomes
\begin{equation}
    \frac{\partial}{\partial \tilde{x}_j} = \frac{1}{s_j(x_j)} \frac{\partial}{\partial x_j}, \quad j = 1, \dots, D.
\end{equation}
Consequently, the Laplacian operator modified by the PML becomes
\begin{equation} 
\tilde{\Delta} = \sum_{j=1}^D \frac{1}{s_j(x_j)} \frac{\partial}{\partial x_j} \left( \frac{1}{s_j(x_j)} \frac{\partial}{\partial x_j} \right) = \sum_{j=1}^D \left( \frac{1}{s_j^2(x_j)} \frac{\partial^2}{\partial x_j^2} - \frac{s_j'(x_j)}{s_j^3(x_j)} \frac{\partial}{\partial x_j} \right). 
\end{equation}
The Helmholtz equation with PML is then written as
\begin{equation}
    -\tilde{\Delta} u(x) - k(x)^2 u(x) = f(x), \quad x \in \Omega \cup \Omega_{\text{abl}}.
\end{equation}

\subsubsection{Fourier Spectral Discretization}
\label{sec:app_fsm}
As an alternative to the finite difference discretization used in the main text (Eq.~\ref{eq:discrete_fdm}), the Fourier spectral method (FSM) discretizes the Helmholtz operator globally by leveraging the spectral representation of the Laplacian:
\begin{equation}
    \Delta u(x)
    = -\mathcal{F}^{-1} \!\left[\, |\xi|^2 \, \mathcal{F}[u(x)] \,\right],
\end{equation}
where $\mathcal{F}$ and $\mathcal{F}^{-1}$ denote the Fourier and inverse Fourier transforms, respectively, and ${\xi}\in \mathbb{R}^D$ is the Fourier wavenumber variable. Applying this to the Helmholtz equation yields a formulation that is exact for the Laplacian term under periodic boundary conditions:
\begin{equation}
    \mathcal{F}^{-1} \left[ |\xi|^2 \, \mathcal{F}[u(x)] \right] - k(x)^2 u(x) = f(x), \quad x \in \mathbb{R}^{D}.
\end{equation}
Discretization via the Fast Fourier Transform (FFT) yields:
\begin{equation}\label{eq:FSM_discrete}
    \mathrm{FFT}^{-1} \left[ |\bm{\xi}|^2 \, \mathrm{FFT}[\bm u] \right] - \bm k^2 \bm u = \bm f.
\end{equation}
Unlike the finite difference method, FSM minimizes dispersion errors by resolving derivatives globally. However, the FFT inherently imposes periodic boundary conditions, which can introduce wrap-around artifacts in non-periodic scattering problems unless combined with sufficient damping (PML).

\subsubsection{Empirical comparison of discretizations and boundaries}
\label{sec:app_discretization_frameworks}
We compare the four combinations of discretization (FDM, FSM) and absorbing boundary (sponge, PML) under the model and dataset configuration of Section~\ref{sec:comp_with_CBS_OpenBreastUS}, with results in Table~\ref{tab:discretization_comp}. Two trends emerge. Within each discretization, the PML degrades both training loss and convergence relative to the sponge layer: it modifies the differential operator itself, replacing the Laplacian by the complex-stretched operator $\tilde\Delta$ derived above, whose additional first-order terms make the learned operator harder to optimize, whereas the sponge layer acts only on the zeroth-order term and leaves the Laplacian unchanged. Across discretizations, FSM is consistently less favorable than FDM under the same boundary, since the periodic boundaries implied by the FFT admit wrap-around error when the attenuation is not sufficiently strong. The two effects compound, making FSM+PML by far the hardest configuration, while FDM with a sponge layer, the setting used throughout the main text, is the most efficient.

\begin{table}[htbp]
\centering
\caption{Effect of discretization (FDM vs.\ FSM) and absorbing boundary (sponge vs.\ PML) on McMg training and inference. Iterations and runtime are averaged over 50 test samples at a relative tolerance of $1\times10^{-6}$.}
\label{tab:discretization_comp}
\begin{tabular}{lccc}
\toprule
Discretization & Train Loss & Avg. Iters & Avg. Runtime (s) \\
\midrule
FDM + Sponge & $1.16\times10^{-2}$ & 7.5  & 0.077 \\
FDM + PML    & $2.72\times10^{-2}$ & 11.3 & 0.135 \\
FSM + Sponge & $3.58\times10^{-2}$ & 10.0 & 0.118 \\
FSM + PML    & $4.96\times10^{-2}$ & 57.6 & 0.647 \\
\bottomrule
\end{tabular}
\end{table}

\subsection{The Born-Series Variant and Its Discrete Green Operator}
\label{sec:app_born_series}
This appendix collects two technical components referenced from the main text: the Born-series formulation of McMg (the Born+McMg variant of Section~\ref{sec:mpmg_iteration_scheme}) and the discrete Green operator used to keep all solvers on a common finite-difference discretization in Section~\ref{sec:comp_with_other_model_OpenBreastUS}.

\subsubsection{McMg for the Born-series system}
For ill-conditioned Helmholtz problems, training McMg directly in the original residual space can lead to unstable optimization and slow convergence. Neural Preconditioned Born Series (NPBS)~\cite{wang2026neural} introduced a unified framework for training neural preconditioners in the preconditioned Born-series residual space, improving both training stability and solver convergence for ill-conditioned Helmholtz systems.

To define this formulation, let $\bm G$ denote the convolution operator with the background Green's function associated with the constant background wavenumber $k_0$, i.e., $\bm G = (-\Delta - k_0^2)^{-1}$, and let $\bm V$ be the scattering potential,
\begin{equation}
    \bm V(\bm x) = k(\bm x)^2 - k_0^2.
\end{equation}
The corresponding Lippmann--Schwinger equation~\cite{lippmann1950variational,chew1999waves} is
\begin{equation}
    \bm u = \bm G(\bm V\bm u+\bm f).
\end{equation}
Equivalently,
\begin{equation}
    (\bm I-\bm G\bm V)\bm u = \bm G\bm f,
    \qquad
    \bm I-\bm G\bm V=\bm G A.
\end{equation}
Iterating this fixed-point equation gives the classical Born iteration. Given an approximate solution $\bm u^{(k)}$, the Born-series residual is
\begin{equation}
    \bm r_{\mathrm{BS}}^{(k)}
    =
    \bm G(\bm V\bm u^{(k)}+\bm f)
    -
    \bm u^{(k)}
    =
    \bm G(\bm f-A\bm u^{(k)}).
\end{equation}
Thus the Born-series residual is precisely the Helmholtz residual after left preconditioning by $\bm G$.
The classical Born iteration may diverge under strong scattering. The Convergent Born Series (CBS)~\cite{osnabrugge2016convergent} restores convergence by introducing a pointwise (diagonal) relaxation operator $\bm\gamma$:
\begin{equation}\label{eq:cbs_eq}
\bm u^{(k+1)}
=
\bm u^{(k)}
+
\bm\gamma
\left(
\bm G(\bm V\bm u^{(k)}+\bm f)
-
\bm u^{(k)}
\right).
\end{equation}
With a suitable choice of $\bm\gamma$, this iteration is guaranteed to converge under the CBS framework, although convergence can still be slow in high-contrast media. Following NPBS, we replace this fixed diagonal relaxation with a learned McMg:
\begin{equation}
    \bm u^{(k+1)}
    =
    \bm u^{(k)}
    +
    \mathcal{MG}
    \left(
    \bm G(\bm V\bm u^{(k)}+\bm f)
    -
    \bm u^{(k)}
    \right).
\end{equation}
The training objective is chosen to match the preconditioned residual space. We sample residuals $\bm r$ as in Section~\ref{sec:mpmg_iteration_scheme} and use $\bm G_{(i)}\bm r$ as the network input, so that training and inference use the same preconditioned coordinates. The corresponding loss is
\begin{equation}
    \mathcal{L}_{\mathrm{BS}}
    =
    \frac{1}{B}
    \sum_{i=1}^{B}
    \frac{
    \left\|
    \left(\bm I-\bm G_{(i)}\bm V_{(i)}\right)
    \mathcal{MG}
    \left(
    \bm k_{(i)},
    \bm G_{(i)}\bm r
    \right)
    -
    \bm G_{(i)}\bm r
    \right\|_2
    }{
    \left\|
    \bm G_{(i)}\bm r
    \right\|_2
    }.
\end{equation}
This trains McMg to approximate the inverse of the \emph{preconditioned} operator $\bm G A$, rather than the inverse of the original Helmholtz operator $A$. The left preconditioning by $\bm G$ regularizes the spectrum of the operator to be inverted, which is what stabilizes training and accelerates convergence in the ill-conditioned regime.

\subsubsection{Discrete Green operator on the finite-difference grid}\label{sec:app_discrete_green}
The Born-series solvers apply $\bm G$ through the FFT, which makes them especially efficient but, in their original form, ties them to the spectral symbol $-|\bm\xi|^2$ of the continuous Laplacian. To compare CBS and Born+McMg against the finite-difference solvers on an identical discrete system (Section~\ref{sec:comp_with_other_model_OpenBreastUS}), we instead build $\bm G$ from the symbol of the discrete Laplacian, retaining the FFT purely as a fast diagonalization tool rather than as a spectral discretization. For the standard five-point stencil in two dimensions with grid spacing $h$, the discrete Laplacian is diagonalized by the Fourier modes with eigenvalues
\begin{equation}\label{eq:discrete_symbol}
\lambda_h(\xi_x,\xi_y)
=
-\frac{4}{h^2}\sin^2\!\left(\frac{\xi_x h}{2}\right)
-\frac{4}{h^2}\sin^2\!\left(\frac{\xi_y h}{2}\right),
\end{equation}
which is the second-order consistent approximation of $-|\bm\xi|^2$ and recovers it as $h\to0$. Replacing $-|\bm\xi|^2$ by $\lambda_h$ in the background Green operator, $\bm G = (-\Delta_h - k_0^2)^{-1}$ is realized as multiplication by $(-\lambda_h - k_0^2)^{-1}$ in Fourier space, so the resulting Green's function is consistent with the five-point finite-difference operator. Consequently, the LU factorization, the shifted-Laplacian preconditioner, and the Born-series methods all act on the same discrete Helmholtz system, and the comparison reflects differences between solvers rather than between discretizations.

\subsubsection{Native continuous-spectrum CBS and Born+McMg}
\label{sec:app_native_cbs}
The main-text comparison (Section~\ref{sec:comp_with_other_model_OpenBreastUS}) deliberately keeps every solver on the common five-point discretization through the discrete Green operator above. For completeness, we also evaluate the Born-series solvers in their \emph{native} form, where $\bm G$ is built from the continuous Laplacian symbol $-|\bm\xi|^2$ and applied via the FFT, rather than from the discrete symbol $\lambda_h$. We retrain Born+McMg in this setting and compare it against native CBS under the OpenBreastUS high-frequency configuration of Section~\ref{sec:comp_with_CBS_OpenBreastUS}, so that the only difference from CBS is the preconditioner. As shown in Table~\ref{tab:native_cbs_comp}, native Born+McMg converges in $6.4$ iterations against $772.1$ for native CBS, a $120.6\times$ reduction in iteration count and a $12.0\times$ wall-clock speedup. These gains match, and slightly exceed, those observed on the common FDM discretization in the main text (Table~\ref{tab:cbs_mpmg_precision}), confirming that the advantage of the learned preconditioner persists in CBS's native continuous-spectrum setting.

\begin{table}[htbp]
\centering
\caption{Native continuous-spectrum Born-series solvers on the OpenBreastUS high-frequency benchmark. Both CBS and Born+McMg build the background Green operator from the continuous Laplacian symbol and apply it via the FFT. Iterations and runtime are averaged over 50 test samples at a relative tolerance of $1\times10^{-6}$; the Ratio column reports CBS relative to Born+McMg.}
\label{tab:native_cbs_comp}
\begin{tabular}{lccc}
\toprule
 & CBS & Born+McMg & Ratio \\
\midrule
Avg. Iters    & 772.1 & 6.4   & 120.6 \\
Avg. Time (s) & 0.311 & 0.026 & 12.0  \\
\bottomrule
\end{tabular}
\end{table}

\subsection{Baseline Solver Configurations}
\label{sec:app_baseline_config}
For reproducibility, we record the parameters of the classical baselines compared in Section~\ref{sec:comp_with_other_model_OpenBreastUS}. All baselines act on the same five-point finite-difference discretization (Appendix~\ref{sec:app_discrete_green}) and run in the same GPU environment; the precision of each method is chosen as described in Section~\ref{sec:comp_with_other_model_OpenBreastUS}, with CBS in FP64 and the remaining iterative methods in mixed precision. Hardware is therefore identical across methods, while the precision of each solver is set to its own fastest configuration so that none is disadvantaged.

\paragraph{Shifted-Laplacian-preconditioned GMRES (GMRES+SL).}
The preconditioner is the complex-shifted Laplacian operator~\cite{erlangga2004class}
\begin{equation}\label{eq:csl_shift}
    M_{(\beta_1,\beta_2)} = -\Delta_h - (\beta_1 + i\,\beta_2)\,k^2,
\end{equation}
with shift $(\beta_1,\beta_2) = (1,\,1.0)$, i.e.\ $M = -\Delta_h - (1+i)\,k^2$. The shifted system is approximately inverted by a single multigrid V-cycle whose smoother is weighted Jacobi with relaxation $\omega = 0.66$ applied to the diagonal of $M$, using $(\nu_1,\nu_2) = (2,2)$ pre- and post-smoothing steps; the coarse level is solved by the same weighted-Jacobi iteration for $2$ steps rather than by a direct solve. The comparatively strong imaginary shift $\beta_2 = 1.0$ (versus the textbook value $0.5$) is adopted precisely because the coarse problem is handled by this inexpensive inner iteration: the added dissipation keeps the V-cycle stable and contractive without an exact coarse solve. GMRES is restarted every $25$ iterations, i.e.\ GMRES($25$), which bounds the Krylov-subspace memory. We note that GMRES+McMg converges in $7.2$ iterations (Table~\ref{tab:cbs_mpmg_precision}), well below the restart length, so restarting is never triggered for McMg and influences only the SL baseline.

\paragraph{Convergent Born Series (CBS).}
CBS is evaluated using the standard aggressive parameterization of Osnabrugge et al.~\cite{osnabrugge2016convergent}. With the scattering potential $\bm V(\bm x) = k(\bm x)^2 - k_0^2$ of Appendix~\ref{sec:app_born_series}, we center the background wavenumber on the real part of the squared-wavenumber range, $k_0^2 = \tfrac12\big(\min_{\bm x} \mathrm{Re}\{k(\bm x)^2\} + \max_{\bm x} \mathrm{Re}\{k(\bm x)^2\}\big)$, the standard choice that reduces the required absorption; for the real-valued $k^2$ of our benchmarks this minimizes $\max_{\bm x}|\bm V(\bm x)|$ exactly. CBS then introduces a uniform absorption $\epsilon$, using the complex potential $\bm V_\epsilon(\bm x) = \bm V(\bm x) - i\epsilon$ and the pointwise diagonal preconditioner
\begin{equation}\label{eq:cbs_gamma}
    \bm\gamma(\bm x) = \frac{i\,\bm V_\epsilon(\bm x)}{\epsilon} = 1 + \frac{i\,\bm V(\bm x)}{\epsilon};
\end{equation}
the same $i\epsilon$ is absorbed into the background Green operator, so the Helmholtz operator $A$ is left unchanged. We take $\epsilon = \max_{\bm x}|\bm V(\bm x)|$, the smallest value permitted by the sufficient CBS convergence condition $\epsilon \ge \max_{\bm x}|k(\bm x)^2 - k_0^2|$ for the no-gain media and absorbing boundaries considered here. The smaller $\epsilon$ is, the longer the effective range of the absorbing background Green's function, and the CBS pseudo-propagation distance per iteration scales as $2k_0/\epsilon$; hence a smaller admissible $\epsilon$ generally leads to faster convergence. This is therefore the most aggressive choice admitted by the CBS formulation, placing CBS in a favorable standard setting rather than at a parameterization proven to be globally optimal.

\subsection{Spectral Analysis of Large Domain Scalability}
\label{sec:app_spectral_analysis}
Consider the constant-coefficient Helmholtz equation $-\Delta u - k^2 u = f$, discretized as the linear system $A \bm u = \bm f$. We aim to learn a parameterized linear operator $\mathcal{MG} \approx A^{-1}$ by minimizing the unsupervised residual objective
\begin{equation}
\min
\; \mathcal{L} = \mathbb E_{\bm r \sim \mathcal N(\bm 0, \bm I)}
\big\| A \mathcal{MG}(\bm r) - \bm r \big\|_2^2.
\label{eq:residual_training_obj}
\end{equation}
To expose the spectral weighting, consider the idealized case in which $A$ is diagonalized by an orthonormal Fourier or sine basis $\{\bm\phi_i\}_{i=1}^N$,
\[
    A\bm\phi_i = \lambda_i \bm\phi_i,
    \qquad
    \lambda_i \simeq |\bm \xi_i|^2-k^2,
\]
where the last relation denotes the corresponding discrete Helmholtz symbol. Boundary and PML effects are suppressed in this local Fourier argument. We further consider the diagonal Fourier response of the learned operator,
\[
    \mathcal{MG}(\bm\phi_i) = m_i \bm\phi_i,
\]
which is exact for a shift-invariant linear operator and serves as the standard local approximation for a trained convolutional multilevel operator. For an input $\bm r=\sum_i \eta_i \bm\phi_i$, with standardized white-noise coefficients satisfying $\mathbb E[\eta_i\overline{\eta_j}]=\delta_{ij}$, the residual error becomes
\begin{equation}
A \mathcal{MG}(\bm r) - \bm r
=
\sum_{i=1}^N \eta_i
\left(
\lambda_i m_i - 1
\right)\bm\phi_i .
\end{equation}
Taking the expectation with respect to $\bm r$ and exploiting the orthogonality of eigenmodes yields:
\begin{equation}
\begin{aligned}
\mathcal{L}
&=
\mathbb E
\left\|
\sum_{i=1}^N
\eta_i
\left(
\lambda_i m_i - 1
\right)\bm\phi_i
\right\|_2^2
\\
&=
\mathbb E \left[
\sum_{i=1}^N \sum_{j=1}^N
\eta_i \overline{\eta_j}
\left(\lambda_i m_i - 1\right)
\overline{\left(\lambda_j m_j - 1\right)}
\left\langle \bm\phi_i, \bm\phi_j \right\rangle
\right]
\\
&=
\sum_{i=1}^N \sum_{j=1}^N
\mathbb E[\eta_i \overline{\eta_j}]
\left(\lambda_i m_i - 1\right)
\overline{\left(\lambda_j m_j - 1\right)}
\left\langle \bm\phi_i, \bm\phi_j \right\rangle
\\
&=
\sum_{i=1}^N
\left| \lambda_i m_i - 1 \right|^2 \\
&=
\sum_{i=1}^N |\lambda_i|^2
\left| m_i - \lambda_i^{-1} \right|^2,
\end{aligned}
\label{eq:spectral_loss}
\end{equation}
for nonzero $\lambda_i$. Thus the residual objective controls the residual error $|\lambda_i m_i-1|$, while the same objective weights the absolute error in the learned inverse response, $|m_i-\lambda_i^{-1}|$, by $|\lambda_i|^2$. For the Helmholtz operator, near-resonant or near-propagating modes satisfy $|\bm \xi_i|^2 \approx k^2$ and hence $|\lambda_i| \ll 1$ in this idealized representation. These modes have large inverse response and are closely associated with the long-range components of the Green's function\cite{sheppard2014green, dwarka2021pollution}. Consequently, finite-capacity and finite-time training can fit the residual objective while still leaving appreciable errors in the inverse action on these near-resonant components. The point is not that the exact minimizer of the residual loss differs from $A^{-1}$; rather, the residual formulation gives weak leverage on the large-gain components that are most important for global wave transport.

On a finite domain of length $L$, the discrete frequency spectrum
has resolution $\Delta \xi = \mathcal{O}(2\pi/L)$. On a small domain, this coarse resolution leaves relatively few near-resonant modes, and the resonance gap $\delta_N=\min_i |\lambda_i|$ is typically larger. Training on such domains is therefore better conditioned and can learn accurate local correction behavior. However, the learned operator is exposed only to wave interactions over the physical scale of the training domain; when transferred to a larger domain, it behaves like a localized or spatially truncated approximation of the Green's function and cannot directly correct all long-range error components. By contrast, as $L$ increases, the spectrum becomes denser near the propagating shell, the number of near-resonant modes grows, and $\delta_N$ decreases. Direct training on the large domain must then learn many large-gain inverse responses that are weakly weighted in the residual-induced inverse-error metric. In practice, this makes optimization less stable and can lead to a learned operator whose residual loss is small but whose global wave coupling remains inaccurate. The resulting long-range error accumulation explains the degradation in convergence observed when the learned correction is used as an iterative solver, and motivates LLPF as a way to expand the learned Green's-function support through newly added coarse levels while preserving the local correction already learned on smaller domains.

\subsection{Ablation Study on Parameter Sharing}
To assess the impact of our architectural design choices on computational efficiency and model capacity, we conduct an ablation study comparing the proposed McMg against two variants with reduced parameter sharing. The variants are defined as follows:
\begin{itemize}
    \item \textbf{McMg:} Adopts the full parameter sharing strategy. The Neural PDE Operator $\mathcal{A}_{l,\bm \theta}$ is shared across all smoothing steps within a level, and physical feature maps ($\bm a_l, \bm s_l$) are pre-computed and cached (Setup Phase).
    \item \textbf{McMg-1:} Disables the sharing of the Neural PDE Operator $\mathcal{A}_{l,\bm \theta}$. Independent operators are instantiated for every smoothing iteration, though the physical feature maps remain shared.
    \item \textbf{McMg-2:} Further disables the sharing of physical feature maps ($\bm a_l, \bm s_l$). This effectively removes the cached Setup Phase, requiring the network to re-compute physical features from the wavenumber $\bm k$ at every iteration.
\end{itemize}

The training characteristics are summarized in Table~\ref{tab:ps_comp_train}. As expected, reducing parameter sharing (McMg-1 and McMg-2) significantly increases the model size, with McMg-2 requiring nearly $3\times$ the parameters of the proposed McMg. While McMg-2 achieves a marginally lower training loss due to its increased expressive capacity, this comes at the cost of a substantially longer training time per epoch (35.22s vs. 21.04s).

Table~\ref{tab:ps_comp_infer} details the inference performance, averaged over 50 test samples with a stopping tolerance of $1 \times 10^{-6}$. The results highlight a critical trade-off: while the highly parameterized McMg-2 requires slightly fewer iterations to converge (7.5 vs. 8.1), its wall-clock inference time is nearly double that of the proposed method (0.091s vs. 0.049s). This confirms that the computational overhead of re-evaluating physical features and applying unshared operators outweighs the marginal gain in convergence speed. The proposed McMg strikes the optimal balance, delivering comparable convergence with significantly higher computational efficiency.

\begin{table}[htbp]
\centering
\caption{Comparison of model complexity and training metrics under different parameter sharing strategies.}
\label{tab:ps_comp_train}
\begin{tabular}{lcccc}
\toprule
Model & Params (MB) & Epoch Time (s) & Train Loss & Val Loss \\
\midrule
McMg   & 2.44 & 21.04 & $1.22 \times 10^{-2}$ & $1.44 \times 10^{-2}$ \\
McMg-1 & 3.81 & 21.06 & $1.24 \times 10^{-2}$ & $1.50 \times 10^{-2}$ \\
McMg-2 & 7.12 & 35.22 & $1.02 \times 10^{-2}$ & $1.25 \times 10^{-2}$ \\
\bottomrule
\end{tabular}
\end{table}

\begin{table}[htbp]
\centering
\caption{Comparison of inference efficiency under different parameter sharing strategies.}
\label{tab:ps_comp_infer}
\begin{tabular}{lcc}
\toprule
Model & Avg. Iterations & Avg. Time (s) \\
\midrule
McMg   & 8.1 & 0.049 \\
McMg-1 & 8.2 & 0.047 \\
McMg-2 & 7.5 & 0.091 \\
\bottomrule
\end{tabular}
\end{table}

\subsection{Ablation study on the number of multigrid levels}
We investigate the effect of varying the number of multigrid levels in the McMg. To ensure a fair comparison, the total number of smoothing steps is fixed across all levels. The detailed configurations are summarized in Table~\ref{tab:exp_levels_train}. For instance, the setting $[1,16]$ corresponds to a two-level model, where level~$1$ performs $1$ smoothing step and level~$2$ performs $16$ steps. To avoid the influence of other factors, the number of channels at each level is fixed to $20$. The wavenumber of the Helmholtz equation and the dataset follow the same settings as in Section~\ref{sec:comp_with_CBS_OpenBreastUS}, while here the discretization is based on finite differences (FDM) for simplicity. As shown in Tables~\ref{tab:exp_levels_train} and \ref{tab:exp_levels_infer}, employing multiple levels achieves better performance under the same total number of smoothing steps, while also improving computational efficiency. Inference results are reported per sample under a relative tolerance of $10^{-6}$.
\begin{table}[htbp]
    \centering
    \caption{Training results for different level configurations}
    \label{tab:exp_levels_train}
    \begin{tabular}{lccc}
        \toprule
        Level Settings & Epoch Time (s) & Train Loss & Val Loss \\
        \midrule
        $[1, 16]$          & 50.38 & $6.23 \times 10^{-2}$ & $6.24 \times 10^{-2}$ \\
        $[1, 4, 12]$       & 43.23 & $4.40 \times 10^{-2}$ & $4.44 \times 10^{-2}$ \\
        $[1, 4, 4, 8]$     & 42.35 & $3.35 \times 10^{-2}$ & $3.49 \times 10^{-2}$ \\
        $[1, 4, 4, 4, 4]$  & 42.56 & $2.76 \times 10^{-2}$ & $2.94 \times 10^{-2}$ \\
        \bottomrule
    \end{tabular}
\end{table}
\begin{table}[htbp]
    \centering
    \caption{Average number of iterations and runtime per sample for different level configurations.}
    \label{tab:exp_levels_infer}
    \begin{tabular}{lcc}
        \toprule
        Level Settings & Avg. Iters & Avg. time(s) \\
        \midrule
        $[1, 16]$          & 366.1 & 2.88 \\
        $[1, 4, 12]$       & 89.6  & 0.58 \\
        $[1, 4, 4, 8]$     & 32.0  & 0.20 \\
        $[1, 4, 4, 4, 4]$  & 21.5  & 0.14 \\
        \bottomrule
    \end{tabular}
\end{table}

\subsection{Sensitivity to Points-Per-Wavelength}
In this experiment, we evaluate the robustness of the proposed McMg architecture under different points-per-wavelength ($\texttt{ppw}$). The difficulty of solving the Helmholtz equation numerically is heavily dependent on the number of $\texttt{ppw}$, we vary the wavenumber to sweep $\texttt{ppw} \in \{6, 8, 10, 12, 14, 16\}$, transitioning from the challenging engineering limit to a well-resolved regime. The model is implemented using FDM discretization and trained on the OpenBreastUS breast dataset, following the protocols in Section~\ref{sec:comp_with_CBS_OpenBreastUS}.

The results, summarized in Table~\ref{tab:exp_diff_wavenumber}, reveal a clear correlation between grid resolution and convergence speed. In the well-resolved regime ($\texttt{ppw} \ge 10$), where the wave is smooth relative to the grid, the McMg converges rapidly, requiring only 4-5 iterations. As $\texttt{ppw}$ decreases towards $6$, the problem difficulty increases sharply due to the pollution effect and numerical dispersion inherent in standard discretizations. While classical geometric multigrid methods typically diverge in this coarse regime due to aliasing of oscillatory modes, McMg maintains stability, converging in 8.2 iterations at $\texttt{ppw} \approx 6$. The increase in iterations at lower $\texttt{ppw}$ indicates that the network must work harder to correct the solver's coarse-grid dispersion mismatch, which sharpens as the wave becomes marginally resolved; this affects only the rate of convergence to the discrete solution, not the accuracy of that solution, and demonstrates the solver's robustness on challenging, marginally resolved cases.

\begin{table}[htbp]
    \centering
    \caption{Performance of McMg across varying $\texttt{ppw}$.}
    \label{tab:exp_diff_wavenumber}
    \begin{tabular}{ccccc}
        \toprule
        $\texttt{ppw}$ & Train Loss & Val Loss & Avg. Iters & Avg. Error \\
        \midrule
        6  & $1.22 \times 10^{-2}$ & $1.44 \times 10^{-2}$ & 8.2 & $4.51 \times 10^{-7}$ \\
        8  & $7.54 \times 10^{-3}$ & $8.98 \times 10^{-3}$ & 5.8 & $3.38 \times 10^{-7}$ \\
        10 & $4.49 \times 10^{-3}$ & $5.44 \times 10^{-3}$ & 5.0 & $1.37 \times 10^{-7}$ \\
        12 & $3.02 \times 10^{-3}$ & $3.53 \times 10^{-3}$ & 4.0 & $2.65 \times 10^{-7}$ \\
        14 & $2.39 \times 10^{-3}$ & $2.79 \times 10^{-3}$ & 4.0 & $1.29 \times 10^{-7}$ \\
        16 & $1.87 \times 10^{-3}$ & $2.16 \times 10^{-3}$ & 4.0 & $5.42 \times 10^{-8}$ \\
        \bottomrule
    \end{tabular}
\end{table}

\end{document}